\documentclass[11pt,a4paper,notitlepage]{article}
\usepackage[]{amssymb}
\usepackage{amssymb}
\usepackage{geometry} 
\usepackage{float}
\usepackage{booktabs}
\usepackage{algorithm,algorithmic}
\usepackage{amsmath}
\usepackage{pdflscape}
\usepackage{adjustbox}
\usepackage{url}
\usepackage{natbib}
\usepackage{graphicx}
\usepackage{lscape}
\usepackage{xcolor}
\usepackage{todonotes}
\usepackage{wasysym} 
\usepackage[utf8]{inputenc}
\usepackage[english]{babel}
\usepackage[normalem]{ulem}
\usepackage{rotating}
\usepackage{setspace}
\usepackage{hyperref}
\usepackage{pifont} 

\usepackage{todonotes}

\allowdisplaybreaks
\title{An exact approach for the multi-depot electric vehicle scheduling problem}
\author{Xenia Haslinger\textsuperscript{a} (xenia.haslinger@jku.at)  \and Elisabeth Gaar\textsuperscript{b} (elisabeth.gaar@uni-a.de) \and Sophie N. Parragh\textsuperscript{a} (sophie.parragh@jku.at) \\ 
\footnotesize \textsuperscript{a} Johannes Kepler University Linz, Institute of Production and Logistics Management, \\
\footnotesize JKU Business School, Altenberger Straße 69, 4040 Linz, Austria \\
\footnotesize \textsuperscript{b} University of Augsburg, Institute of Mathematics, \\
\footnotesize Universitätsstraße 14, 86159 Augsburg, Germany }
\date{}

\begin{document}

\maketitle

\begin{abstract}
    The ``avoid - shift - improve'' framework and the European Clean Vehicles Directive set the path for improving the efficiency and ultimately decarbonizing the transport sector. While electric buses have already been adopted in several cities, regional bus lines may pose additional challenges due to the potentially longer distances they have to travel. 
    
    In this work, we model and solve the electric bus scheduling problem, lexicographically minimizing the size of the bus fleet, the number of charging stops, and the total energy consumed, to provide decision support for bus operators planning to replace their diesel-powered fleet with zero emission vehicles. We propose a graph representation which allows partial charging without explicitly relying on time variables and derive 3-index and 2-index mixed-integer linear programming formulations for the multi-depot electric vehicle scheduling problem. While the 3-index model can be solved by an off-the-shelf solver directly, the 2-index model relies on an exponential number of constraints to ensure the correct depot pairing. These are separated in a cutting plane fashion. 

    We propose a set of instances with up to 80 service trips to compare the two approaches, showing that, with a small number of depots, the compact 3-index model performs very well. However, as the number of depots increases, the developed branch-and-cut algorithm proves to be of value. The inclusion of realistic instances and technology-specific scenarios for diesel, battery-electric, and fuel cell-electric buses further strengthens the practical relevance of our results, offering concrete guidance for sustainable fleet planning.\\
    
\end{abstract}

\section*{Acknowledgments}
This project is supported with funds from the Climate and Energy Fund and
implemented in the framework of the RTI-initiative “Flagship region Energy”, FFG Project
Nr. FO999900979, associated to the energy
model region WIVA P\&G.

\section*{Declaration of competing interest}
The authors declare that they have no known competing financial interests or personal relationships that could have appeared to influence the work reported in this paper.

\newpage

\section{Introduction}

CO$_2$-emissions from transport accounted for about 21\% of the global CO$_2$-emissions in 2023,  increasing by almost 80\% between 1990 and 2023 \citep{statista2024carbon,statista2024transport}.
Following the ``avoid - shift - improve'' framework to meet the Paris Climate Agreement, Austria's 2030 Mobility Masterplan promotes active modes of transport, like walking and biking, public modes of transport, and shared mobility as well as decarbonizing the transportation sector in general \citep{masterplan}. The implementation of the European Clean Vehicles Directive \citep{cleanVehicles} requires all member states to meet at least half of their procurement targets for clean buses via zero emission buses. E.g., Austria's target for heavy duty vehicles (including buses) is currently 45\% and will rise to 65\% with January 2026. 
The Austrian research project ``Zero Emission Mobility Salzburg'' (ZEMoS) is a collaborative effort of several public and private institutions and research organizations to support the transition from diesel-powered to zero emission public bus transport and waste collection vehicles. The main focus is the decarbonization of the regional public bus system  \citep{zemos,innov2024,jahrbuchKapitel}. 

In order to support the transformation of the public bus system in European countries and beyond, high-quality decision support is necessary, especially, concerning the bus operators' important strategic decision of how many and which types of zero emission buses to invest in. Given an available line plan, switching from diesel-powered buses to zero emission vehicles without changes in the plan, requires solving an (electric) vehicle scheduling problem, minimizing the number of buses required to serve the line. Since battery-electric zero emission vehicles are currently still only available with rather limited driving ranges compared to diesel-powered vehicles, and have the disadvantage of taking a long time to recharge, it is not straightforward to substitute conventional diesel buses with zero emission technology, in particular for regional bus lines. 
Recharging stops during the day are usually required and have to be planned.
From an operator's perspective, charging on their own premises, potentially relying on their own photovoltaic energy, is in many cases the cheapest and sometimes the only option, since chargers for heavy duty vehicles are not yet available on a larger scale, resulting in a rather small number of possible charging locations.

The problem we address can be cast as a Multi-Depot Electric Vehicle Scheduling Problem (MDEVSP): each bus line operator has to serve a set of service trips with a fleet of electric or fuel cell-electric vehicles. Each service trip has to be covered by exactly one schedule, which is performed by one vehicle. Each vehicle is required to begin and end its trip at the same depot, chosen from a predefined set of depots. 
Given the use of electric vehicles and the regional setting with relatively long bus routes, charging may need to be planned throughout the day to ensure that the battery level never drops below a specified threshold and remains above a certain limit upon returning to the depot.
Recharging requires vehicles to drive to a charging station, where charging facilities are usually available. Partial recharging is allowed. 

The aim of our approach is to provide decision support for transitioning from diesel-powered to zero emission vehicles, particularly electric buses. The key issue to resolve is determining how many buses are necessary to serve all service trips.
Service trips follow a timetable and, therefore, they have fixed start and end times, durations, 
and lengths. We also assume that the fleet is limited to a single vehicle type, i.e., buses of a specific length and technology. This reflects the situation currently encountered in Salzburg, where individual lines are tendered and assigned to different operators based on their bids. 
However, various bus types can be evaluated and compared to select the most suitable option. Since minimizing the total number of vehicles usually results in several alternative optimal solutions, we
lexicographically minimize (1) the number of vehicles, (2) the number of charging events during the day, and (3) the total energy consumed, while all scheduled service trips during daily operations are covered.

The contribution of this paper is as follows:
\begin{itemize}
\item We develop a graph representation that only allows time-feasible paths (or vehicle schedules), including charging events.
\item  We formulate a 3-index and a 2-index mixed-integer linear program (MILP), relying on our graph representation, allowing partial recharge without explicitly modeling time.
\item The 3-index model is compact and can be solved directly with an off-the-shelf solver. The 2-index model requires depot-pairing constraints of exponential size. We propose two different ways to model these constraints (infeasible path constraints and connectivity constraints) and separate them in a branch-and-cut fashion.
\item We present a computational study on a large set of instances that mimic realistic settings and illustrate the advantages and drawbacks of the two approaches.
\item We compare propulsion technologies (diesel, battery-electric, and fuel cell-electric) using realistic instances, conducting a parameter study to examine the impact of cold weather conditions as well as battery degradation on fleet size and operational feasibility. 
\end{itemize}

The remainder of this paper is structured as follows. We first give an overview of the work related to the problem setting considered in this paper in Section~\ref{sec:lit}. We then formally define the MDEVSP with partial recharging and detail a graph representation of it in Section~\ref{sec:problemFormulation}. Next, we present two mathematical formulations for the MDEVSP, one with three and one with two indices in Section \ref{sec:milps}. We propose a branch-and-cut algorithm to solve this problem in Section~\ref{sec:branchandcut}. The algorithm is tested on newly generated MDEVSP benchmark instances based on an instance generation scheme commonly used in the literature. Computational results are provided in Section \ref{compResults}, followed by 
a comparison of propulsion technologies in Section \ref{techComparison}.
Conclusions are drawn in Section \ref{conclusion}.

\section{Literature review}
\label{sec:lit}

Electric vehicle scheduling problems have attracted quite some attention over the past years. In the following, we review some of the most closely related contributions, starting with the work of \cite{Li_2014} who studied the single-depot electric vehicle scheduling problem (EVSP) without partial charging, and have shown that it is NP-hard. They develop a branch-and-price approach, considering battery swapping or fast charging up to battery capacity in conjunction with maximum distance constraints. Capacitated charging stations in conjunction with partial recharging in the context of electric vehicle scheduling in public transit is considered by \citet{de2024electric}. They rely on a path-based formulation which is solved by column generation. Integer solutions are obtained by price-and-branch and a diving heuristic. Instances with up to 816 service trips are addressed.

\citet{Janovec2019} develop a MILP for the single-depot problem with electric vehicles, considering possible charging events at each of the available chargers. Their formulation requires 4-index variables and, similar to our approach, they do not rely on time variables explicitly. The model is applied to case study-based data with up to 160 service trips and solved by Xpress IVE.

\citet{Wen2016} present a MILP for the MDEVSP considering partial charging relying on an almost acyclic network representation but requiring time variables, minimizing a combination of travel and vehicle costs. Larger instances are addressed by a large neighborhood search algorithm. 
\citet{Wang2021} develop a column generation algorithm and combine it with a genetic algorithm to address the MDEVSP of three bus lines in Qingdao, China.
\citet{Adler2017} address the the MDEVSP with partial recharge, minimizing the total schedule costs, limiting the number of buses,  which can be stationed at each depot. A branch-and-price algorithm is developed and used to benchmark a heuristic approach on random instances with up to 50 service trips. The heuristic is then applied to a large-scale data set from the metropolitan area of Phoenix.

\citet{Hu_2022} address the optimization of locating fast chargers, which allow en-route charging at selected bus stops along three bus routes in Sydney, and determining charging schedules assuming time-dependent electricity prices. Further they consider the possibility to delay service, but penalize passenger waiting times. Passenger demand as well as travel times are assumed to be uncertain. The developed robust optimization model is solved by Gurobi. Time-of-use electricity prices are also considered by \citet{vanKooten2017scheduling}, \citet{Li2020}, and \citet{wu2022multi} in a deterministic setting. \cite{vanKooten2017scheduling} present two mathematical programming models for the single-depot case, where in the first, a linear charging process with constant electricity prices during the day is assumed. The second model allows any type of charging process, includes time-of-use electricity prices, and takes the depreciation cost of the battery into account. For the latter model, the exact value charged is approximated by discretizing it. Instances with up to 175 service trips are solved to optimality with both models, for larger instances with 241 service trips, methods based on column generation are deployed. \citet{wu2022multi} formulate the problem with two objectives, which are considered in a lexicographic fashion and develop a branch-and-price algorithm. 
Also \citet{liu2020battery} address the combined problem of vehicle scheduling and charger location, allowing partial recharge for transit buses in Shanghai. The developed modeling approaches rely on deficit function theory and mixed-integer programming. 
The location of charging infrastructure in combination with electric bus scheduling and a single depot is modeled as a 2-index MILP by \citet{Stumpe2024} and, using Gurobi as the solver engine, applied to large-scale instances.

\cite{zhou2024electric} propose a MILP and a set-covering formulation for the electric bus charging scheduling problem with a mixed fleet of electric vehicles, partial recharging while taking battery degradation effects into account. A branch-and-price algorithm solves instances with up to 100 service trips, for tackling large-scale instances, an optimization-based adaptive large neighborhood search (ALNS) is developed. 
A mixed fleet of conventional and electric vehicles is also considered by \citet{Sassi2017}. Since charging costs change over time, in the developed model, time is discretized and costs are minimized in the objective function. Also, maximum grid capacity constraints are considered.
\citet{Rinaldi2020} address mixed fleet bus scheduling with a single depot. Time is discretized, assuming that buses can be fully charged within one time interval. Trip departure times may deviate from their preferred time. The developed MILP is decomposed into smaller subproblems and applied to case study data from Luxembourg City.
\citet{Yildirim2021} also address a mixed fleet problem, but with multiple depots and with multiple different charging technologies and develop an efficient column generation-based algorithm to determine the minimum cost fleet configuration for large-scale real-world instances. 
\cite{zhang2022mixed} present a mixed-integer programming formulation for the MDEVSP with heterogeneous vehicles. They solve the model using CPLEX and propose an ALNS algorithm to address the problem, incorporating a partial mixed-route strategy and a partial recharging policy.
\citet{friess2024planning} have developed a MILP, where different zero emission propulsion technologies are considered concurrently to optimize the bus fleet mix for serving the city of Graz, Austria. In order to solve the model, they pre-select a set of possible options for buses to switch between different lines.

\citet{Zhang2021b} address electric bus fleet scheduling with a single depot, non-linear charging and battery fading, and develop a branch-and-price approach. Non-linear charging, considering battery degradation effects are also studied by \citet{Zhou2022}.
A column-generation based heuristic enhanced by ideas from machine learning has recently been developed by \citet{Gerbeau2025} and applied to the MDEVSP with non-linear charging. Also \citet{Diefenbach2022} consider non-linear charging in an in-plant vehicle scheduling application of the multi-depot problem. In their case, vehicles are free to return to any depot, not necessarily the one they started from. They address this problem setting by a branch-and-check approach, moving all complicating aspects, such as the planning of charging events, to the subproblem and generate cutting planes, enforcing a change in the current solution, when violated. 
\cite{lobel2024electric} propose a MILP for the electric bus scheduling problem with a single depot for mixed fleets of electric and non-electric vehicles, presenting an improved approximation of non-linear battery charging behavior of the electric vehicles. Furthermore, they address the challenge of adjusting to power grid bottlenecks by integrating dynamic recharging rates and time-of-use electricity prices. Partial charging is allowed and available charging slots must not be exceeded. The model is applied to diverse real-life instances with up to 1,207 service trips.

\cite{jiang2022branch} present a mixed-integer program and develop a branch-and-price algorithm to address real-world instances of the MDEVSP with up to 460 service trips, considering a partial recharging policy, time windows for service trip start times, and charging depots located in close proximity to the start and end stops of each line. The proposed branch-and-price algorithm is improved by incorporating a heuristic approach to generate good initial solutions as well as by embedding heuristic decision-making within the label-setting algorithm to solve the pricing problem. Also the branching step is enhanced by heuristic rules for fixing variables to 0 or 1, if their fractional values are very close to these values. Single line instances without time windows and up to 200 service trips are solved to proven optimality.

\cite{Gkiotsalitis2023} propose a mixed-integer non-linear model, which they later reformulate to a MILP, for the electric MDVSP with time windows, where operational cost of buses as well as vehicle waiting time are considered and service trips may start within a certain time window.  
In their formulation, simultaneously charging different vehicles on the same charger is prohibited. The authors introduce valid inequalities to tighten the search space of the MILP. The implementation is demonstrated on a toy network and on randomly generated benchmark instances similar to \cite{Carpaneto1989}. 
For further work on electric bus scheduling we refer to the survey conducted by \cite{Perumal2022}.

\section{Graph representation of the MDEVSP}
\label{sec:problemFormulation}

On the way to our ultimate goal of solving the MDEVSP by utilizing mixed-integer linear programming, we first lay the foundations in this section. In particular, we formally introduce the MDEVSP in Section~\ref{sec:definitionMDEVSP}, and then we detail in Section~\ref{sec:graphRepresentation} how any instance of the MDEVSP can be transformed into a graph, such that solving a particular kind of flow problem in this graph is equivalent to solving the MDEVSP for this instance.

\subsection{Description of the MDEVSP}
\label{sec:definitionMDEVSP}
We start by giving the formal definition of the MDEVSP.
In the MDEVSP, we are given a set of timetabled service trips $V^I$, where each service trip $i \in V^I$ corresponds to a scheduled trip with a specific start time~$s_i$ and end time~$e_i$, a duration~$u_i$ and an energy usage~$q_i$. Moreover, each service trip $i$ has a start location $\ell^s_i$ and an end location~$\ell^e_i$.

Furthermore, we are given the set of depot indices $K$ and for each depot index $k \in K$ we are given the number of available vehicles $b_k$ and the depot location 
$\ell_{k}$. 
Note that the depots for different depot indices can be at the same physical location.
For notational convenience we introduce 
the set of origin depots $O$ consisting of $o_k$ for all $k \in K$, and the set of destination depots $D$ consisting of $d_k$ for all $k \in K$. 
For each depot index $k \in K$, we define the start and the end location of the origin and destination depot as $\ell_{k}$, so 
$\ell^s_{o_k} = \ell^e_{o_k} = \ell^s_{d_k} = \ell^e_{d_k} = \ell_{k}$ hold.

Moreover, we are given the set of charging stations $C$, where each charging station $a \in C$ has a location $\ell_a$, which serves as both start and end location, so $ \ell^s_a = \ell^e_a = \ell_a$.
The locations of the charging stations can be at any specified physical location and may coincide with the depot locations. 

For each pair of service trips, origin or destination depots, or charging stations $i,j \in V^I \cup O \cup D \cup C$, we are given 
the traveling time from  $\ell^e_i$ to $\ell^s_j$ (i.e., the time it takes a vehicle to travel from the end location of $i$ to the start location of $j$) as $t_{ij}$, 
and the distance between  $\ell^e_i$ and $\ell^s_j$ as $d_{ij}$.
Moreover, we compute the energy usage $p_{ij}$ between the locations $\ell^e_i$ and $\ell^s_j$ as $p_{ij}= \theta d_{ij}$, where $\theta$ is the given energy consumption rate, i.e., the amount of energy used per distance unit. Thus, a vehicle's energy consumption is assumed to be a linear function of the traveled distance.

Furthermore, the maximum battery capacity of the electric vehicles $s^{max}$, the minimum allowed battery level $s^{min}$, the minimum allowed battery level when returning to the depot at the end of the schedule $s^{min}_{dep}$ (with $s^{min}_{dep} \geq s^{min}$), the minimum time that needs to be available for charging $t^{min}$ and the charging rate $r$, i.e., the time needed to charge one unit of the battery, are given. From that we can compute the time required to fully charge the battery $t^{max}$ as $t^{max} = (s^{max}-s^{min})/r$.
A vehicle's charging time is assumed to be a linear function of the amount that its battery is charged.

The goal of the MDEVSP is to find a schedule, such that 
(1) each timetabled service trip $i \in V^I$ is done by exactly one vehicle, (2) each vehicle starts in some depot $k \in K$, performs a sequence of service trips and returns to the same depot $k$ at the end of the day, 
(3) between two service trips or between a service trip and returning to the depot, each vehicle visiting a charging station can reload its battery for at least $t^{min}$ minutes with charging rate $r$, 
(4) each vehicle starts at its depot with energy level $s^{max}$ and returns to its depot with energy level at least $s^{min}_{dep}$ at the end of the day, 
(5) the energy level of each vehicle never drops below $s^{min}$, and
(6) in each depot $k \in K$ at most $b_k$ vehicles start and end their routes. 
We call such a schedule feasible.
Among all feasible schedules the MDEVSP searches for the one which lexicographically minimizes (i) the number of used vehicles, (ii) the number of visits of charging stations, and (iii) the energy spent on deadhead trips. The latter is equivalent to minimizing the total energy consumed, as the energy spent on service trips is a constant value.

\citet{bertossi1987some} have shown that the vehicle scheduling problem with a single depot can be solved in polynomial time while its multi-depot version, the MDVSP, is NP-hard. Since the MDVSP is a special case of the MDEVSP, where the battery capacity of each vehicle $s^{max} = \infty$ and there are no charging stations, the MDEVSP is NP-hard as well.

Each feasible schedule represents the planned itinerary for all vehicles throughout the day.
The schedule of one vehicle can be represented as an ordered sequence of service trips, depots and charging stations, so as an ordered sequence of elements from $V^I \cup O \cup D \cup C$. 
In particular, the first element of the sequence must belong to the set of origin depots $O$, the last to the set of destination depots $D$, and all others in between to $V^I \cup C$. Additionally, the sequence must include at least one service trip, no two consecutive charging stations can occur and each charging station can occur only directly after a service trip. 
This illustration of schedules as ordered sequences will be the foundation of our graph representation of the MDEVSP.

\subsection{Graph representation}
\label{sec:graphRepresentation}

Our next step is to define a graph, 
such that each flow in this graph with certain properties corresponds to a feasible solution of the MDEVSP.
To construct the graph, 
initially we set the set of nodes $V$ as $V = V^I \cup V^D \cup V^C$, where the service trip node set $V^I$, the depot node set $V^D = O \cup D$, and the charging node set $V^C := \emptyset$ are used. Furthermore, we start with the arc set $A := \emptyset$. 
We add arcs to $A$ and charging nodes to $V^C$ by following these steps:

\begin{enumerate}
    \item Generate arcs between depots and service trips:
    For each depot index $k \in K$, we generate an arc from the origin depot $o_k$ to each service trip $i \in V^I$, as well as from each service trip $i \in V^I$ to the destination depot $d_k$. 

    \item Generate an arc between each pair of time-feasible service trips: 
    For each pair of service trips $i,j \in V^I$, we create an arc if $s_i + u_i + t_{ij} \leq s_j$ holds. We store these pairs of time-feasible service trips in $F$.
    
    \item Generate (and partially connect) full charging nodes representing fully charging after service trips: 
    For each service trip $i \in V^I$ and each charging station $a \in C$ we create a charging node $c^{full}_{ia} \in V^C$ representing the possibility of fully charging at charging station $a$ after service trip $i$. 
    This charging node $c = c^{full}_{ia}$ has location $\ell_c = \ell_a$,
    maximum available charging time $t_{c} = t^{max}$ and maximum amount that can be charged $h_{c} = r t_c$.
    Additionally, we connect $i$ and $c$ with an arc and this arc has energy usage $p_{ic} = p_{ia}$.

    \item Generate charging nodes and connect (full) charging nodes representing charging between two service trips: 
    For each pair of service trips $(i,j) \in F$, we identify 
    the set $C_{ij} \subseteq C$ of all charging stations, which are reasonable for charging a vehicle between the service trips $i$ and $j$. To be more precise, a charging station $a \in C$ is in $C_{ij}$, if there is no other charging station $a' \in C$ that dominates $a$. A charging station $a'$ dominates a charging station $a$, if (1) $a'$ is closer to the end location $\ell^e_i$ of service trip $i$ than $a$ and (2) $a'$ is closer to the start location $\ell^s_j$ of service trip $j$ than $a$ (if such a charging station $a'$ exists, then it would always be better to charge at $a'$ than at $a$).
    
    For each charging station $a \in C_{ij}$, 
    we compute the maximum available charging time $t_{ija} = \min \{s_j - (s_i + u_i + t_{ia} + t_{aj}), t^{max} \}$ and the maximum amount that can be charged $h_{ija} = r t_{ija}$.
    If the energy required to travel from service trip $i$ to service trip $j$ via charging station $a$ is less than the amount that can be recharged at charging station $a$, so if $p_{ia} + p_{aj} < h_{ija}$, 
    we do the following case distinction.

    On the one hand, if we can fully charge, so if $t_{ija} = t^{max}$, 
    then we use the full charging node $c = c^{full}_{ia}$.
    We connect $c$ and $j$ with an arc and we set $p_{cj} = p_{aj}$.

    On the other hand, if we can not fully charge, but we can charge at least $t^{min}$ minutes, so if $t^{max} > t_{ija} \geq t^{min}$, then
    we generate a charging node $c^{part}_{ija} \in V^C$ representing charging at $a$ between the service trips $i$ and $j$. 
    Charging node $c=c^{part}_{ija}$ has location $\ell_c = \ell_a$,
    maximum available charging time $t_{c} = t_{ija}$ and maximum amount that can be charged $h_{c} = h_{ija}$.
    Additionally, we connect both $i$ and $c$ and $c$ and $j$ with an arc and set the energy usage $p_{ic} = p_{ia}$ and $p_{cj} = p_{aj}$.

    \item Generate arcs representing fully charging before returning to the depot: 
    For each service trip $i \in V^I$ and each destination depot $d_k \in D$, we determine
    the set $C_{i{d_k}} \subseteq C$ of all charging stations, which are reasonable for charging a vehicle between service trip $i$ and going back to the depot $d_k$. In particular, a charging station $a \in C$ is in $C_{i{d_k}}$, if there is no other charging station $a' \in C$ such that (1) $a'$ is closer to the end location $\ell^e_i$ of service trip $i$ than $a$ and (2) $a'$ is closer to the depot location $\ell_{d_k}$ of depot $d_k$ than $a$ (because if such a charging station $a'$ exists, then it would always be better to charge at $a'$ than at~$a$).
    For each charging station $a \in C_{i{d_k}}$, 
    we then connect the full charging node $c= c_{ia}^{full} \in V^C$ 
    and $d_k$ with an arc and set the energy usage $p_{c,{d_k}} = p_{a,{d_k}}$.   
\end{enumerate}

This gives us the graph $G=(V,A)$.
For illustration purposes, we provide a small example instance $I_1$. We consider three service trips $V^{I} = \{ST_1, ST_2, ST_3\}$, one depot (resulting in the start depot $o_1$ and the end depot $d_1$), and two charging stations $C = \{a_1,a_2\}$. Furthermore, we assume that the vehicles start fully charged from the depot.

\begin{figure}[htb]
\centering
\includegraphics[width=0.7\textwidth]{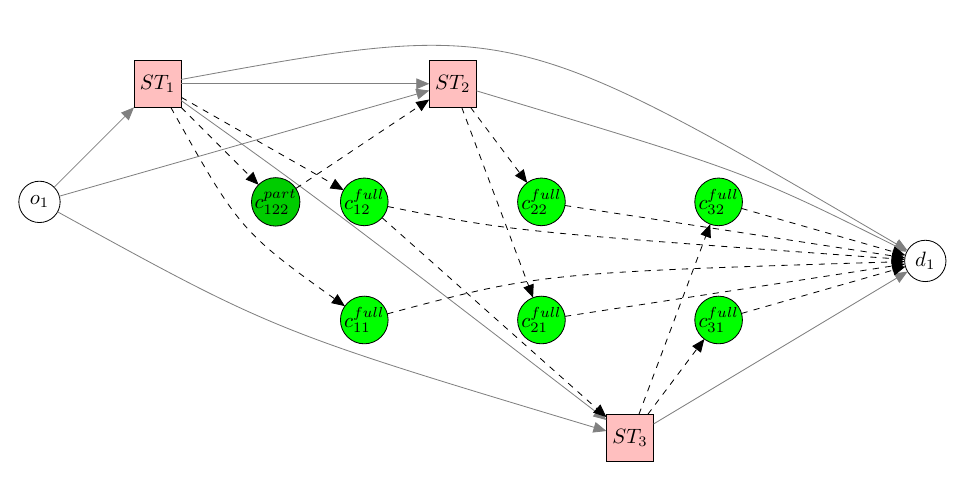} 
\caption{Graph for the instance $I_1$ with a single origin ($o_1$) and destination ($d_1$) depot.}
\label{fig:SD-graph}
\end{figure}

The travel time $t_{ij}$ between service trips, depots, and charging stations, as well as the start time $s_i$ and end time $e_i$ of each service trip can be found in Table~\ref{tab:t_ijGraph}. The time $t^{max}$ it takes to fully charge a vehicle at a full charging node $c^{full}_{ia}$ is set to $120$ minutes. Recharging at the partial charging nodes $c^{part}_{121}$ and $c^{part}_{122}$ is possible for $72$ and $116$ minutes, respectively. However, since charging station $a_2$ is closer to both $ST_1$ and $ST_2$ than $a_1$, $ST_1$ and $ST_2$ are only connected via $c^{part}_{122}$. The same holds for $ST_1$ and $ST_3$, where enough time is available to fully recharge, but again $a_2$ is closer to both $ST_1$ and $ST_3$ than $a_1$. Therefore, $ST_1$ and $ST_3$ are only connected via $c^{full}_{12}$. An example of the graph for instance $I_1$ is provided in Figure~\ref{fig:SD-graph}.

\begin{table}[h!]
    \centering
        \caption{Traveling time $t_{ij}$ between service trips, depots, and charging stations, as well as the service trips' start and end times $s_i$ and $e_i$ for the example graph in Figures~\ref{fig:SD-graph} and \ref{fig:MD-graph}.}
    \label{tab:t_ijGraph}
    \begin{tabular}[t]{c| c c c|c c|c c }
    \toprule
        $t_{ij}$ &  $ST_1$ & $ST_2$ & $ST_3$ & $o_1/d_1$ & $o_2/d_2$  & $a_1$ & $a_2$ \\ \hline
        $ST_1$ & - & 28 & 5 & 34 & 15 & 28 & 19 \\ 
        $ST_2$ & 28 & - & 5 & 34 & 15 & 28 & 19  \\ 
        $ST_3$ & 35 & 35 & - & 7 & 49 & 32 & 23 \\ \hline
        $o_1/d_1$ & 40 & 40 & 29 & - & 47 & 26 & 26 \\ 
        $o_2/d_2$ & 39 & 39 & 19 & 47 & - & 33 & 34 \\ \hline
        $a_1$ & 50 & 50 & 24 & 20 & 30 & - & 36 \\
        $a_2$ & 15 & 15 & 16 & 26 & 34 & 35 & - \\ 
    \end{tabular}
    \quad
    \begin{tabular}[t]{c| c | c }
    \toprule
        ~ &  $s_i$ & $e_i$ \\ \hline
        $ST_1$ & 01:15 p.m. & 02:00 p.m.\\ 
        $ST_2$ & 04:30 p.m. & 05:15 p.m.\\ 
        $ST_3$ & 05:05 p.m. & 06:30 p.m. \\ 
    \end{tabular}
\end{table}

Grey arcs are direct connections between service trips and between service trips and depots. Black dashed arcs connect service trips with charging nodes, as well as charging nodes with the depot. All connections are time feasible. E.g., $ST_2$ and $ST_3$ are not connected via an arc, since $ST_2$ ends at 5:15 p.m. but $ST_3$ already starts at 5:05 p.m. All service trips are connected via two full charging nodes with the destination depot $d_1$, since charging station $a_1$ is closer to the depot while $a_2$ is closer to the service trips. Thus, connections via full charging nodes have been generated for all combinations of service trips and charging stations.

\begin{figure}[htb]
\centering
\includegraphics[width=0.7\textwidth]{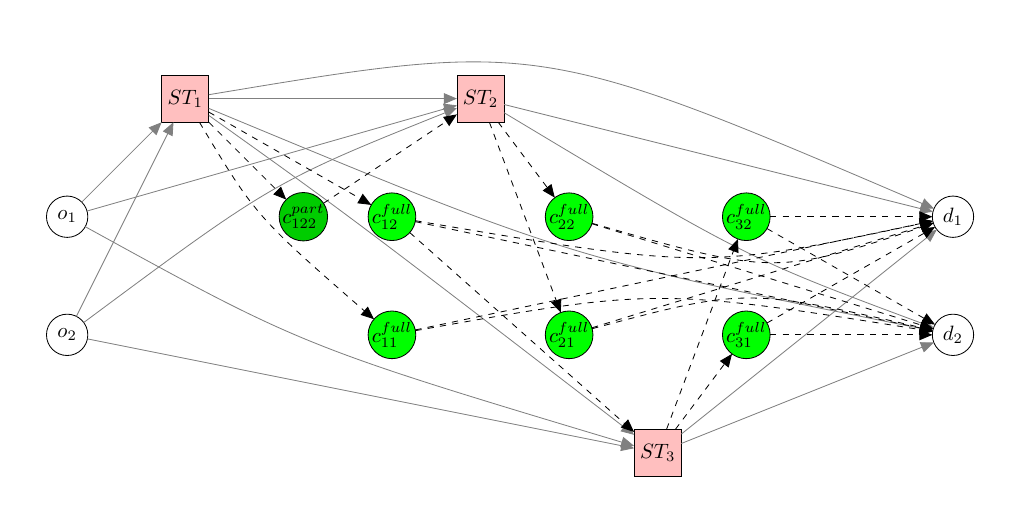} 
\caption{Graph for the instance $I_2$ with two origin depots ($o_1, o_2$) and two destination depots ($d_1, d_2$).}
\label{fig:MD-graph}
\end{figure}

For another example instance $I_2$, which is identical to the instance $I_1$, except for the addition of an extra depot, we examine a small example of the corresponding graph provided in Figure~\ref{fig:MD-graph}. The service trips as well as the reasonable full charging nodes are now connected to both depots. The rest of the graph structure remains unchanged from Figure~\ref{fig:SD-graph}.

Eventually, we have defined the graph $G=(V,A)$ such that each feasible solution of the MDEVSP corresponds to a certain type of flow (namely one that makes sure that each unit of flow ends at the same depot as it starts) in $G$. In particular, in such a flow each unit of flow represents the schedule of one vehicle and starts from an origin depot $o_k \in O$ for some $k \in K$ and ends at the corresponding destination depot $d_k \in D$.
Note that we have embedded all time relevant information into the graph construction directly, 
such that any flow in our graph will automatically correspond to a time feasible schedule for the MDEVSP.

\section{Mathematical models for the MDEVSP}
\label{sec:milps}

The graph representation of the MDEVSP will now be the foundation of two mixed-integer linear programming formulations of the MDEVSP. First, we present a 3-index formulation in Section~\ref{sec:3indexFormulation}. Then we detail how we can omit one of the indices and derive a 2-index formulation in Section~\ref{sec:2indexFormulation}. Furthermore, we present valid inequalities and possible extensions in Section~\ref{sec:validinequalities}.
Finally, we describe how our formulations can be adapted to alternative zero emission technologies in Section~\ref{sec:extensions}.

\subsection{A 3-index formulation}
\label{sec:3indexFormulation}

We now present our first mixed-integer linear programming formulation for the MDEVSP based on the graph $G=(V,A)$ derived in Section~\ref{sec:graphRepresentation}, which we refer to as 3-index formulation. Our models are inspired by \citet{Stumpe2021} and \citet{Friess2021}. 

For notational convenience, we will denote by $A^-(i)$ and $A^+(i)$ the set of all predecessor and successor nodes of a node $i \in V$ in $G$, respectively. So
$A^-(i) = \{j \in V: (j,i) \in A\}$ and $A^+(i) = \{j \in V: (i,j) \in A\}$ holds.

We introduce two sets of decision variables for this 3-index formulation: For each $k\in K$ and each $(i,j) \in A$ let $x_{ij}^k \in \{0,1\}$ be a binary variable, which is 1, if the arc $(i,j)$ is used from vehicles associated with depot $k$, 0 otherwise. Furthermore, for each $i \in V$, let $\varepsilon_i$ be the energy level of vehicles when leaving node $i$.
With these decision variables, the 3-index MDEVSP can be formulated as the following mixed-integer multi-commodity flow model.
\begin{subequations}
\label{z3_MDEVSP_model}
\begin{align}
    \min~ w_1 \Big(\sum_{k \in K} \sum_{j \in A^+(o_k)} & x^k_{o_k,j}\Big) 
    + w_2 \Big(\sum_{k \in K} \sum_{c \in V^C} \sum_{i \in A^-(c)}  x^k_{ic} \Big)
    + && w_3 \Big(\sum_{k \in K} \sum_{(i,j) \in A} p_{ij} x^k_{ij}\Big) 
    \label{z3_MDEVSP_of}
    \\
\text{s.t.~} 
    \sum_{k \in K} \sum_{i \in A^-(j)} x_{ij}^k & = 1 && \forall j \in V^I &
    \label{constr:z3_MDEVSP_visit}\\
    \sum_{j \in A^+(o_k)} x^k_{o_k,j} &\leq b_k  && \forall k \in K &
    \label{constr:z3_MDEVSP_buslimit}\\
    \sum_{i \in A^-(j)} x^k_{i j} & =  \sum_{i \in A^+(j)} x^k_{j i} && \forall j \in V \setminus V^D, k \in K &&
    \label{constr:z3_MDEVSP_inout}\\
    \sum_{j \in A^+(o_k)} x^k_{o_k,j} & = \sum_{i \in A^-(d_k)} x^k_{i,d_k}  && \forall k \in K &
    \label{constr:z3_MDEVSP_samenumberbuses}\\
    \sum_{\substack{k' \in K\\ k' \neq k}} \sum_{j \in A^+(o_{k'})} 
    x^k_{o_{k'},j} & = 0 && \forall k \in K  &
    \label{constr:z3_MDEVSP_samebusoutfromO}\\
    \sum_{\substack{k' \in K\\ k' \neq k}} \sum_{i \in A^-(d_{k'})}     
    x^k_{i, d_{k'}} & = 0 && \forall k \in K &
    \label{constr:z3_MDEVSP_samebusinD}\\
    \varepsilon_{o_k} & = s^{max} && \forall k \in K &
    \label{constr:z3_MDEVSP_energyOrigin}\\
        s^{min} & \leq \varepsilon_i && \forall i \in V \setminus V^D &
    \label{constr:z3_MDEVSP_minEnergy}\\
        s^{min}_{dep} & \leq \varepsilon_{d_k} &&  \forall k \in K & 
        \label{constr:z3_MDEVSP_minEnergyDepot}  \\  
    s^{min} & \leq  \varepsilon_i - \sum_{j \in A^+(i)} p_{ij}x^k_{ij}  &&  \forall i \in V^I, k \in K
    \label{constr:z3_MDEVSP_minEnergyDT}\\
    \varepsilon_j & \leq \varepsilon_i -  (p_{ij}+q_j) x^k_{ij} + s^{max}(1-x^k_{ij})  && \forall j \in V^I, i \in A^-(j), k \in K &
    \label{constr:z3_MDEVSP_consumeEnergyST}\\
    \varepsilon_c & \leq \varepsilon_i -(p_{ic} - h_c)  x^k_{ic} + s^{max} (1-x^k_{ic}) && \forall c \in V^C, i \in A^-(c), k \in K
     \label{constr:z3_MDEVSP_loadEnergy}\\
    \varepsilon_{d_k} & \leq \varepsilon_i - p_{i,d_k} x^k_{i,d_k}  + s^{max} (1-x^k_{i,d_k}) && \forall k \in K, i \in A^-(d_k)
    \label{constr:z3_MDEVSP_energyDestination}\\
    \varepsilon_c & \leq s^{max} && \forall c \in V^{C} 
    \label{constr:z3_MDEVSP_maxEnergy} \\
    x_{ij}^k & \in \{0,1\} && \forall (i,j) \in A, k \in K
    \label{constr:z3_MDEVSP_domainX} 
\end{align}
\end{subequations}

The objective function~\eqref{z3_MDEVSP_of} minimizes, in a lexicographic way, first the number of vehicles (which equals the number of used arcs leaving all origin depots $o_k$), second the number of charging events during the day (which equals the sum of the used arcs entering any of the charging nodes $c$), and third the energy spent on deadhead trips (which equals the required energy for all used arcs). 
Towards this end, the objective function uses the weights $w_1$, $w_2$, $w_3 \in \mathbb{R}$, which must be set appropriately and fulfill $w_1 > w_2 > w_3$.

Constraints~\eqref{constr:z3_MDEVSP_visit} make sure that each service trip node $j \in V^I$ is visited exactly once.
Constraints~\eqref{constr:z3_MDEVSP_buslimit} limit the number of vehicles that can be used at each depot.
Constraints~\eqref{constr:z3_MDEVSP_inout} ensure flow conservation (if a node is entered it has to be left unless it is a depot).
Constraints~\eqref{constr:z3_MDEVSP_samenumberbuses} make sure that at each depot the same number of vehicles leaves and arrives.
Constraints~\eqref{constr:z3_MDEVSP_samebusoutfromO} guarantee that vehicles associated with the depot with index $k$ use only arcs leaving origin depot $o_k$ and no arcs leaving another origin depot $o_{k'}$ for some $k' \neq k$, while constraints \eqref{constr:z3_MDEVSP_samebusinD} ensure that vehicles associated with the depot with index $k$ use only arcs arriving at destination depot $d_k$ and no arcs arriving at another destination depot $d_{k'}$ for some $k' \neq k$.

Constraints~\eqref{constr:z3_MDEVSP_energyOrigin} set the energy level at all origin depots to the maximum energy level (all vehicles leave the depots fully charged).
Constraints~\eqref{constr:z3_MDEVSP_minEnergy} make sure that the energy level when leaving any node except the depots is not below the minimum required energy level.
Constraints~\eqref{constr:z3_MDEVSP_minEnergyDepot} ensure that the energy level when returning to the depot is not below the required energy level at the destination depot.
Constraints~\eqref{constr:z3_MDEVSP_minEnergyDT} make sure that the energy level is not below the minimum required energy level before arriving at the charging station (which is implied by the fact that the energy level is not below the minimum required energy level when arriving at the next node after a service trip node, as charging nodes are always preceded by service trip nodes).
Constraints~\eqref{constr:z3_MDEVSP_consumeEnergyST} guarantee that the energy level at the end of a service trip corresponds to the energy level of the previous node minus the energy consumed by the deadhead trip connecting the previous node and this service trip and the energy used by the service trip itself. 
Constraints~\eqref{constr:z3_MDEVSP_loadEnergy} allow the energy level at the end of a charging node to rise to at most the energy level of the previous node minus 
the energy that was consumed by the deadhead trip connecting the previous node and the charging station plus the energy that can be loaded at this charging node.
Constraints~\eqref{constr:z3_MDEVSP_energyDestination} make sure that the energy level at the destination depots corresponds to at most the energy level at the previous node minus the energy spent on the deadhead trip connecting the previous node to the depot.
Constraints~\eqref{constr:z3_MDEVSP_maxEnergy} make sure that a vehicle cannot be charged to an energy level that is higher than the battery capacity. 
The domain of the variable $x$ is defined in~\eqref{constr:z3_MDEVSP_domainX}.

Note that \eqref{z3_MDEVSP_model} is a compact mixed-integer linear program, that can be solved with off-the-shelf solvers. We will present computational results for doing so later on.

\subsection{A 2-index formulation}
\label{sec:2indexFormulation}

Next, we present a second mixed-integer linear programming formulation for the MDEVSP. Like our first formulation \eqref{z3_MDEVSP_model}, it is based on the graph $G=(V,A)$ derived in Section~\ref{sec:graphRepresentation}. 

For this new formulation, which we refer to as 2-index formulation, we introduce a binary decision variable $x_{ij} \in \{0,1\}$ for each $(i,j) \in A$. It is 1 if arc $(i,j)$ is used by any vehicle, 0 otherwise. The second set of decision variables $\varepsilon_i$ is used in the same sense as in the 3-index formulation, i.e., it is the energy level of vehicles when leaving node $i$ for each $i \in V$.
Then the 2-index formulation of the MDEVSP is given as  
\begin{subequations}
\label{z2_MDEVSP_model}
\begin{align}
    \min~ w_1 \Big(\sum_{k \in K} \sum_{j \in A^+(o_k)} & x_{o_k,j}\Big) 
    + w_2 \Big(\sum_{c \in V^C} \sum_{i \in A^-(c)}  x_{ic} \Big)
    + && w_3  \Big(\sum_{(i,j) \in A}   p_{ij} x_{ij} \Big) 
    \label{z2_MDEVSP_of}
    \\
\text{s.t.~} 
    \sum_{i \in A^-(j)} x_{ij} & = 1 && \forall j \in V^I &
    \label{constr:z2_MDEVSP_visit}\\
    \sum_{j \in A^+(o_k)} x_{o_k,j} &\leq b_k  &&  \forall k \in K&
    \label{constr:z2_MDEVSP_buslimit}\\
    \sum_{i \in A^-(j)} x_{i j} & =  \sum_{i \in A^+(j)} x_{j i} && \forall j \in V \setminus V^D &&
    \label{constr:z2_MDEVSP_inout}\\
    \sum_{j \in A^+(o_k)} x_{o_k,j} & = \sum_{i \in A^-(d_k)} x_{i,d_k}  && \forall k \in K&
    \label{constr:z2_MDEVSP_samenumberbuses}\\
        \sum_{(i,j) \in P} x_{ij} & \leq |P| - 1  && \forall P \in \mathcal{P}
    \label{constr:z2_MDEVSP_infeasiblePath} \\
    \varepsilon_{o_k} & = s^{max} &&  \forall k \in K &
    \label{constr:z2_MDEVSP_energyOrigin}\\
            s^{min} & \leq \varepsilon_i && \forall i \in V \setminus V^D &
    \label{constr:z2_MDEVSP_minEnergy}\\
        s^{min}_{dep} & \leq 
        \varepsilon_{d_k} &&  \forall k \in K & 
        \label{constr:z2_MDEVSP_minEnergyDepot}  \\  
    s^{min} & \leq  \varepsilon_i - \sum_{j \in A^+(i)} p_{ij}x_{ij}  &&  \forall i \in V^I
    \label{constr:z2_MDEVSP_minEnergyDT}\\
    \varepsilon_j & \leq \varepsilon_i - (p_{ij}+q_j) x_{ij}  + s^{max}(1-x_{ij})  && \forall j \in V^I, i \in A^-(j) &
    \label{constr:z2_MDEVSP_consumeEnergyST}\\
    \varepsilon_c & \leq \varepsilon_i - (p_{ic} - h_c) x_{ic} + s^{max} (1-x_{ic})  && \forall c \in V^C, i \in A^-(c)
     \label{constr:z2_MDEVSP_loadEnergy}\\
    \varepsilon_{d_k} & \leq \varepsilon_i - p_{i,d_k} x_{i,d_k}  + s^{max} (1-x_{i,d_k}) && \forall  k \in K, i \in A^-(d_k)
    \label{constr:z2_MDEVSP_energyDestination}\\
    \varepsilon_c & \leq s^{max} && \forall c \in V^{C} 
    \label{constr:z2_MDEVSP_maxEnergy} \\
    x_{ij} & \in \{0,1\} && \forall (i,j) \in A,
    \label{constr:z2_MDEVSP_domainX}
\end{align}
\end{subequations}
where $\mathcal{P}$ is the set of all paths $P = (v_1, ..., v_n)$ in $G$ with $v_1$, $\ldots$, $v_n \in V$ and $(v_i, v_{i+1}) \in A$ for all $i = 1, \ldots, n-1$ from an origin depot $v_1 = o_k$ for some $k \in K$ to another other destination depot $v_n = d_{k'}$ for some $k' \in K\setminus \{k\}$. Thus, the so-called infeasible path constraints~\eqref{constr:z2_MDEVSP_infeasiblePath} \citep[see, e.g.,][]{Ascheuer2000polyhedral} ensure that no vehicle takes an infeasible path (starting and ending in a different depot). 

All other constraints of the 2-index formulation~\eqref{z2_MDEVSP_model} are defined analogously to the constraints of the 3-index formulation~\eqref{z3_MDEVSP_model}, with the difference that now we do not have a separate variable $x^k_{ij}$ for each $k \in K$ representing vehicles associated to the depot $k \in K$, but only one variable $x_{ij}$ for each arc $(i,j) \in A$.
Thus, it is not possible to directly tell from the variables of~\eqref{z2_MDEVSP_model} at which depot a vehicle using the arc $(i,j)$ departed.
This advantage of having fewer variables comes at the cost of having a high number of infeasible path constraints~\eqref{constr:z2_MDEVSP_infeasiblePath}. We describe in Section~\ref{sec:separationInfeasiblePath} how to deal with this computationally.

Another option to avoid the undesirable infeasible paths is to replace constraints \eqref{constr:z2_MDEVSP_infeasiblePath} with other constraints similar to how it is done, e.g., in \cite{Parragh2011}. 
For doing so, for each $k \in K$ we introduce the set $\mathcal{U}_k$ of all node subsets $U \subseteq V$, such that the origin depot $o_k \in U$ and the destination depot $d_k \notin U$, while for all other $k' \in K\setminus\{k\}$ the origin depots $o_{k'} \notin U$ and the destination depots $d_{k'} \in U$.
Then, we can replace the infeasible path constraints \eqref{constr:z2_MDEVSP_infeasiblePath} in the 2-index formulation~\eqref{z2_MDEVSP_model} with 
\begin{equation}
      \sum_{i \in U}~ \sum_{j \in A^+(i), j \notin U} x_{ij} \geq \sum_{j \in A^+(o_k)} x_{o_k,j} \qquad \forall k \in K, \, U \in \mathcal{U}_k.
    \label{val_i:z2_MDEVSP_maxFlow}   
\end{equation}
These connectivity constraints~\eqref{val_i:z2_MDEVSP_maxFlow} ensure that for every $k \in K$ every vehicle leaving the origin depot $o_k$ must exit the set $U$ (which contains $o_k$) and enter the set $V \setminus U$ (which contains~$d_k$), guaranteeing the correct pairing of the corresponding depots. 
Note that the cardinality of each of the sets $\mathcal{U}_k$ is $2^{|V^I \cup V^C|}$, as $V = V^I \cup V^D \cup V^C$, so there are exponentially many connectivity constraints~\eqref{val_i:z2_MDEVSP_maxFlow}.
We detail in Section~\ref{sec:separationConnectivityCon} how we treat them in our computations.

\subsection{Valid inequalities}
\label{sec:validinequalities}

We now turn our attention to improving the MILP models~\eqref{z3_MDEVSP_model} and~\eqref{z2_MDEVSP_model} by adding valid inequalities.
Our first valid inequalities are based on a lower bound of vehicles needed for covering all the service trips. In particular, we determine the maximum number of concurrent service trips $LB$, which is defined as the maximum number of service trips that are timetabled at the same time. This count $LB$ provides us with a lower bound on the required number of vehicles, i.e., the number of vehicles leaving any of the origin depots $o_k$. Thus, 
\begin{subequations}
\begin{align}
     \sum_{k \in K} \sum_{j \in A^+(o_k)} x_{o_k,j}^k  \geq LB
    \label{val_i:z3_MDEVSP_lowerboundbuses}  \\
\intertext{is a valid inequality for the 3-index formulation \eqref{z3_MDEVSP_model} and}
     \sum_{k \in K} \sum_{j \in A^+(o_k)} x_{o_k,j}  \geq LB
    \label{val_i:z2_MDEVSP_lowerboundbuses}  
\end{align}
\end{subequations}
is a valid inequality for the 2-index formulation \eqref{z2_MDEVSP_model}.

Our second set of valid inequalities for the 3-index formulation~\eqref{z3_MDEVSP_model} is based on decreasing the constant of the big-$M$-type constraints \eqref{constr:z3_MDEVSP_consumeEnergyST}, \eqref{constr:z3_MDEVSP_loadEnergy} and \eqref{constr:z3_MDEVSP_energyDestination} that ensure that the energy level of each vehicle is propagated through the graph in the right way. In particular, $\varepsilon_i \geq s^{min}$ is always fulfilled for each $i \in V$ because of \eqref{constr:z3_MDEVSP_energyOrigin}, \eqref{constr:z3_MDEVSP_minEnergy} and \eqref{constr:z3_MDEVSP_minEnergyDepot}. Thus, the constant $s^{max}$ in the constraints \eqref{constr:z3_MDEVSP_consumeEnergyST}, \eqref{constr:z3_MDEVSP_loadEnergy} and \eqref{constr:z3_MDEVSP_energyDestination} can be replaced by the smaller constant $(s^{max} - s^{min})$ and hence \eqref{constr:z3_MDEVSP_consumeEnergyST}, \eqref{constr:z3_MDEVSP_loadEnergy} and \eqref{constr:z3_MDEVSP_energyDestination} can be strengthened to 
\begin{subequations}
\begin{align}
\varepsilon_j & \leq \varepsilon_i - (p_{ij}+q_j) x_{ij}^k  + (s^{max}-s^{min})(1-x_{ij}^k)  && \forall j \in V^I, i \in A^-(j),  k \in K \label{val_i:z3_MDEVSP_consumeEnergyST} \\
    \varepsilon_c & \leq \varepsilon_i - (p_{ic} - h_c) x_{ic}^k + (s^{max}-s^{min})(1-x_{ic}^k)  && \forall c \in V^C, i \in A^-(c),  k \in K
     \label{val_i:z3_MDEVSP_loadEnergy}\\
 \varepsilon_{d_k} & \leq \varepsilon_i - p_{i,d_k} x_{i,d_k}^k  + (s^{max}-s^{min})(1-x_{i,d_k}^k) && \forall k \in K, i \in A^-(d_k).
    \label{val_i:z3_MDEVSP_energyDestination}
\end{align}
\end{subequations}
With this smaller constant the same integer solutions remain feasible for the 3-index formulation \eqref{z3_MDEVSP_model}, while the feasible region of the linear relaxation becomes smaller, leading to hopefully stronger LP bounds.

Analogously, one obtains valid inequalities for the 2-index formulation~\eqref{z2_MDEVSP_model} by  
\begin{subequations}
\begin{align}
\varepsilon_j & \leq \varepsilon_i - (p_{ij}+q_j) x_{ij}  + (s^{max}-s^{min})(1-x_{ij})  && \forall j \in V^I, i \in A^-(j)  \label{val_i:z2_MDEVSP_consumeEnergyST} \\
    \varepsilon_c & \leq \varepsilon_i - (p_{ic} - h_c) x_{ic} + (s^{max}-s^{min})(1-x_{ic})  && \forall c \in V^C, i \in A^-(c)
     \label{val_i:z2_MDEVSP_loadEnergy}\\
 \varepsilon_{d_k} & \leq \varepsilon_i - p_{i,d_k} x_{i,d_k}  + (s^{max}-s^{min})(1-x_{i,d_k}) && \forall k \in K, i \in A^-(d_k).
    \label{val_i:z2_MDEVSP_energyDestination}
\end{align}
\end{subequations}

\subsection{Extension to alternative zero emission technologies}
\label{sec:extensions}

Finally, we want to point out that even though we have created the graph and derived the two MILP formulations~\eqref{z3_MDEVSP_model} and~\eqref{z2_MDEVSP_model} for the MDEVSP considering electric charging at dedicated charging stations, it is also possible to apply our framework for alternative zero emission technologies.

For example, our model is able to depict the possibility of opportunity charging at the start or end locations of service trips. To do so, for each service trip $i \in V^I$, if opportunity charging is available at the start (end) location of the service trip $\ell_i^s$ ($\ell_i^e$), then a charging station $a$ is added to $C$ at location $\ell_a = \ell_i^s$ ($\ell_a = \ell_i^e$) and the corresponding values of $t_{ij}$, $d_{ij}$ and $p_{ij}$ for $i=a$ or $j=a$ need to be adapted accordingly. Only then the graph is constructed.

Moreover, the option of overnight charging at a depot $k \in K$ can be integrated into the model by creating a charging station $a$ at the location $\ell_k$ of the depot $k$.

Furthermore, fuel cell-electric vehicles can be considered with our methodology if the charging stations $a\in C$ represent hydrogen fueling stations, and the energy consumption rate $\theta$, the charging rate $r$, the maximum energy level $s^{max}$, the minimum energy level $s^{min}$ and $s^{min}_{dep}$ and the minimum time that needs to be available for charging $t^{min}$ are modified to fit for the fuel cell case.

Thus, our approach is universal in the sense that it can be adapted to various zero emission technology settings. Also, the consideration of diesel vehicles is possible with our model analogously to the fuel cell case.

\section{Branch-and-cut algorithm}
\label{sec:branchandcut}

In order to employ the previously introduced 3-index formulation~\eqref{z3_MDEVSP_model} and 2-index formulation~\eqref{z2_MDEVSP_model} to solve the MDEVSP, one could use one of the many MILP solvers available in the standard configuration. While this is possible for~\eqref{z3_MDEVSP_model}, the large number of infeasible path constraints~\eqref{constr:z2_MDEVSP_infeasiblePath} or connectivity constraints \eqref{val_i:z2_MDEVSP_maxFlow} becomes prohibitive for~\eqref{z2_MDEVSP_model}.

Thus, we have developed a branch-and-cut algorithm for~\eqref{z2_MDEVSP_model}. Branch-and-cut algorithms incorporate the principles of branch-and-bound and pair it with the cutting-plane idea. They start from solving the linear relaxation of the MILP, while considering only a reasonable subset of the original constraints. Typically, constraints of exponential size are excluded. Then, in the course of the algorithm, a separation method is required, which finds violated constraints and adds them in an iterative fashion, until ultimately no original constraints are violated anymore, even though they might not be included explicitly.

In our branch-and-cut algorithm, we start with solving~\eqref{z2_MDEVSP_model} with all constraints except for the infeasible path constraints \eqref{constr:z2_MDEVSP_infeasiblePath}, which are initially omitted. 
We call this model the base model from now on. 
We now explore two options for the branch-and-cut algorithm: either we separate and add infeasible path constraints~\eqref{constr:z2_MDEVSP_infeasiblePath}, or connectivity constraints \eqref{val_i:z2_MDEVSP_maxFlow} in the course of our algorithm.

\subsection{Separation of infeasible path constraints}
\label{sec:separationInfeasiblePath}

We start by investigating the separation of infeasible path constraints~\eqref{constr:z2_MDEVSP_infeasiblePath}. Whenever we are given a feasible solution $(x,\varepsilon)$ to the base model (which implies integrality of the variables $x_{ij}$ for each arc $(i,j) \in A$ because of~\eqref{constr:z2_MDEVSP_domainX}), we check if there is a path $P=(v_1, \dots, v_n)$ in $G=(V,A)$ with $v_1, \dots, v_n \in V$ and $(v_i,v_{i+1}) \in A$ such that $x_{v_i,v_{i+1}} = 1$ for all $i = 1, \dots, n-1$ that starts in a depot $v_1 = o_k$ for some $k \in K$ and ends in another depot $v_n = d_{k'}$ for some $k' \in K \setminus \{k\}$. Such a path $P$ corresponds to a vehicle arriving at a different depot than it started and hence is infeasible, which implies that $(x,\varepsilon)$ is feasible for the base model, but infeasible for the 2-index formulation~\eqref{z2_MDEVSP_model}. As a result, such an infeasible path must be prohibited. 

Clearly, whenever one infeasible path is found, at least one other path is violated, as the number of vehicles departing and arriving at each depot is the same. Thus, we consider two different options: either we add only one  infeasible path constraint~\eqref{constr:z2_MDEVSP_infeasiblePath} as soon as we find the first infeasible path $P$ (option \texttt{One}), or we collect all infeasible paths and add the constraint~\eqref{constr:z2_MDEVSP_infeasiblePath} for all of them (option \texttt{All}). 

In the so far described separation we only separate whenever an integer solution is encountered within the branch-and-cut algorithm. We will refer to this setting as \texttt{I}. It is possible to additionally use a separation in the case a fractional solution (i.e, a feasible solution to the linear relaxation of the current node problem in the branch-and-bound tree) is encountered. 
In particular, whenever we are given a feasible solution $(x,\varepsilon)$ to the linear relaxation of the base model, we 
check if there is a path $P=(v_1, \dots, v_n)$ in $G=(V,A)$ with $v_1, \dots, v_n \in V$ and $(v_i,v_{i+1}) \in A$ such that $x_{v_i,v_{i+1}} > 0.00001$ for all $i = 1, \dots, n-1$ that starts in a depot $v_1 = o_k$ for some $k \in K$ and ends in another depot $v_n = d_{k'}$ for some $k' \in K \setminus \{k\}$, and such that the infeasible path constraint~\eqref{constr:z2_MDEVSP_infeasiblePath} is violated for $P$.
Whenever such an infeasible path is found, we add it in the same fashion as in setting \texttt{I}.
We refer to this setting of separating for both integer and fractional solutions as \texttt{IF}. 
Note, that this kind of separation is optional in the sense that even though no fractional solutions are separated, the branch-and-cut algorithm is still correct, as only integer feasible solutions need to be separated for correctness.

\subsection{Separation of connectivity constraints}
\label{sec:separationConnectivityCon}

Next, we draw our attention to separating the connectivity constraints \eqref{val_i:z2_MDEVSP_maxFlow}.
For a feasible solution $(x,\varepsilon)$ to the base model (where $x_{ij}$ is binary for each arc $(i,j) \in A$), we first check for a depot with index $k \in K$, whether there is an infeasible path starting in depot $o_k$ and ending in a different depot $d_{k'}$ for some $k' \in K\setminus \{k\}$ like it is done in the separation for infeasible path constraints described in Section~\ref{sec:separationInfeasiblePath}.
If we have identified a depot with index $k \in K$ where an infeasible path starts, we want to identify a set $U \in \mathcal{U}_k$ for which the connectivity constraint~\eqref{val_i:z2_MDEVSP_maxFlow} is violated by a large amount. To do so, we solve a max flow problem with source $o_k$ and sink $d_k$ in the graph $G=(V,A)$, where the capacity of an arc $(i,j)$ is exactly the value of $x_{ij}$. This allows to determine a minimum cut $(T,V\setminus T)$ with $T\subseteq V$, $o_k \in T$ and $d_k \in V\setminus T$ by the max-flow min-cut theorem such that 
\begin{equation}
    \sum_{i \in T}~ \sum_{j \in A^+(i), j \notin T} x_{ij}
\end{equation}
is minimized, and thus $U = T$ can be used for adding a connectivity constraint~\eqref{val_i:z2_MDEVSP_maxFlow}.
Note, that some technical modification is necessary for $G$ in order to make sure that $o_{k'} \in V\setminus T$ and $d_{k'} \in T$ for all $k' \in K \setminus\{k\}$.

Again, due to the fact that if some path is infeasible, there are at least two depots $k\in K$ such that an infeasible path starts in $k$, we have the option of stopping the separation and adding the constraint as soon as one infeasible path is found (option \texttt{One}), or we can determine a set $U$ for all depots $k\in K$ in which infeasible paths start and add a constraint~\eqref{val_i:z2_MDEVSP_maxFlow} for all such depots (option \texttt{All}).

Additionally, we consider both the classical setting \texttt{I} of separating only integer solutions, and separating both integer and fractional solutions \texttt{IF} analogously as in the separation of infeasible path constraints. Again, both versions ensure correctness of our branch-and-cut algorithm.

\section{Computational experiments}
\label{compResults}

We are now able to present computational results. In the following, first, the attributes of the generated test instances are described in Section~\ref{sec:comp:instances}. Then, our obtained results are discussed for instances with a single depot in Section~\ref{sec:comp:one}, a low number of depots (two, three, and four) in Section~\ref{sec:comp:low}, and a high number of depots (six and eight) in Section~\ref{sec:comp:high}. 

\subsection{Computational setup and test instances}
\label{sec:comp:instances}
Everything is implemented in Julia 1.11.1. CPLEX 22.1 is employed for the branch-and-cut algorithms, where for the separation as described in Section~\ref{sec:branchandcut} we utilize \texttt{LazyConstraint} for separating integer solutions and 
\texttt{UserCut} for separating fractional solutions.
In the separation of the connectivity constraints~\eqref{val_i:z2_MDEVSP_maxFlow} described in Section~\ref{sec:separationConnectivityCon}, we use the Boykov-Kolmogorov algorithm within the function \texttt{maximum\_flow()} of Julia in order to determine a minimum cut.
All experiments were carried out on a Quad-core X5570 Xeon CPU @2.93GHz with a memory of 48 GB.

For our computational experiments, we use a set of MDEVSP benchmark instances generated in a similar way to the generation of MDVSP class A instances in \cite{Carpaneto1989}, used by many other researchers (see \cite{Gkiotsalitis2023,Bianco1994,Fischetti2001,Forbes1994,Ribeiro1994}). The instances of \cite{Carpaneto1989} are not directly applicable to our setting, as they do not involve electric vehicles, which necessitate the inclusion of charging infrastructure in the instance data. Thus, we have created instances as close as possible to~\cite{Carpaneto1989} and adapted them to the electric case, as it was done in~\cite{Gkiotsalitis2023}.

In particular, 
we consider instances 
with $|V^I| \in \{10, 20, 30, 40, 50, 60, 70, 80\}$ service trips, $|K| = |O| = |D| \in \{1,2,3,4,6,8\}$ depots, and $|C| \in \{1,2,3, 4, 6\}$ charging stations.
For each instance, we first determine a number $\nu$ of potential start and end locations of service trips (so-called relief locations) by choosing $\nu$ as uniformly random integer in $\left[ |V^I|/3, |V^I|/2 \right]$. 
Then we choose the $\nu$ relief locations $\ell^{ST}_1$, \dots, $\ell^{ST}_\nu$,
the depot locations $\ell_k$ for each $k\in K$, as well as the charging station locations $\ell_a$ for each $a \in C$ randomly distributed in a 60 km by 60 km square in the Euclidean plane using a uniform distribution, resulting in coordinates (latitude and longitude) for each location. 
Finally, for each service trip $i \in V^I$, we decide if it is a short trip (with a probability of $40\%$, representing urban journeys) or a long trip (with probability $60\%$, representing extra-urban journeys that start and end at the same location). Then for each service trip $i \in V^I$ we choose both the start and the end location of the service trip $\ell^s_i$ and $\ell^e_i$ uniformly at random as one of the $\nu$ relief locations $\ell^{ST}_1$, \dots, $\ell^{ST}_\nu$, where we make sure that $\ell^s_i = \ell^e_i$ holds for long service trips $i$.
From these locations, we compute the distances $d_{ij}$ as the Euclidean distances between the (end) location of $i$ and the (start) location of $j$ for each $i,j \in V^I \cup O \cup D \cup C$.
For the travel times $t_{ij}$ (in minutes), we assume an average vehicle speed of 60 km per hour, translating to 1 km per minute, so $t_{ij} = d_{ij}$.

For each short service trip $i\in V^I$, we generate the start time $s_i$ (in minutes since midnight) as a random integer in the interval [420,480] with a probability of 15\%, in the interval [480,1020] with a probability of 70\% and in the interval [1020,1080] with a probability of 15\%. Moreover, we choose the end time $e_i$ uniformly at random as integer in the interval $[s_i+d_i+5, s_i+d_i+40]$, where $d_i$ is the Euclidean distance between $\ell^s_i$ and $\ell^e_i$.
For each long service trip $i\in V^I$ we generate the start time $s_i$ as a uniform random integer in the interval [300,1200] and the end time $e_i$ as a uniform random integer in the interval  $[s_i+180, s_i+300]$.
For each service trip $i \in V^I$, this yields the duration $u_i$, from which we compute the energy usage $q_i$ as $q_i=\theta u_i$.

For each depot index $k \in K$, we generate the number of vehicles available as $b_k$ as a uniformly random integer in $[3 + \frac{1}{3|K|} V^I,
3 + \frac{1}{2|K|}|V^I|]$.
The parameters related to vehicles and charging are set to $\theta= 1.3$, $s^{max}=1000$, $s^{min}=10$, and $r=50/6$, as done in~\cite{Gkiotsalitis2023}. Note, that this yields $t^{max}=(s^{max}-s^{min})/r = 118.8$. 
Furthermore, we use $s^{min}_{dep}=0.7s^{max}$ and $t^{min}=s^{max}/100$.

For each instance, we first generate the graph $G=(V,A)$ as described in Section~\ref{sec:graphRepresentation}. As this can be done quite fast (39.63 seconds for the largest instance) we omit the running time for doing so from now on. This yields graphs with a broad range of sizes, ranging from 
$|V| = 29$ and $|A| = 172$ for one of the smallest instances with $|V^I|=10$, $|K|=1$ and $|C|=1$ 
up to
$|V| = 2,187$ and $|A| = 22,624$ for one of the largest instances with $|V^I|=80$, $|K|=8$ and $|C|=6$. Note that the majority of the nodes of $G$ is charging nodes (for example, $|V^C| = 2,091$ for the latter instance with $|V| = 2,187$).

The graphs $G=(V,A)$ enable us to utilize the 3-index formulation~\eqref{z3_MDEVSP_model} and the 2-index formulation~\eqref{z2_MDEVSP_model} and its variant.
We carefully engineered $w_1 = 100,000$, $w_2 = 4,000$ and $w_3 = 1$ as appropriate weights for  the weights of the objective functions of~\eqref{z3_MDEVSP_model} and~\eqref{z2_MDEVSP_model} for all our instances, to make sure that we indeed minimize the desired quantities (number of vehicles, number of charging events, energy usage for deadhead trips) in lexicographic order.

\subsection{Results for a single depot}
\label{sec:comp:one}
We start by investigating instances with only one depot, so $|K| = |O| = |D| = 1$. For these instances the 3-index formulation~\eqref{z3_MDEVSP_model} and the 2-index formulation~\eqref{z2_MDEVSP_model} coincide, as neither infeasible path constraints~\eqref{constr:z2_MDEVSP_infeasiblePath} (the set $\mathcal{P}$ is empty), nor connectivity constraints \eqref{val_i:z2_MDEVSP_maxFlow} (the set $\mathcal{U}_k$ is empty) are present.

We consider a set of 90 instances, namely five instances for each combination of $|V^I| \in \{10, 20, 30, 40, 50, 60\}$ service trips and $|C| \in \{1,2,3\}$ charging stations.

\begin{table}[!htb]
    \centering
    \small
    \caption{Results for a single depot.}
    \label{tab:3index_oneDepot}
    \begin{tabular}{rr|rrrrrr}
    \toprule
$|V^I|$& $|C|$ & \#opt & $t$(s) & gap(\%)  & \#nB\&B & $z^*$ & $LB_R$ \\ \midrule
10 & 1 & 5 & 0.03 & 0.00 & 0.00 & 473,330.79 & 472,278.81 \\ 
10 & 2 & 5 & 0.05 & 0.00 & 0.00 & 475,884.69 & 475,190.80 \\ 
10 & 3 & 5 & 0.05 & 0.00 & 0.00 & 555,755.03 & 555,511.01 \\ \hline
20 & 1 & 5 & 0.04 & 0.00 & 2.20 & 866,918.71 & 866,913.57 \\ 
20 & 2 & 5 & 0.31 & 0.00 & 631.60 & 850,045.99 & 845,011.42 \\ 
20 & 3 & 5 & 0.29 & 0.00 & 310.80 & 765,020.56 & 761,797.69 \\  \hline
30 & 1 & 5 & 1.33 & 0.00 & 2,217.80 & 1,117,568.26 & 1,104,410.40 \\ 
30 & 2 & 5 & 0.42 & 0.00 & 413.40 & 1,158,258.26 & 1,150,531.03 \\ 
30 & 3 & 5 & 2.80 & 0.00 & 2,182.80 & 1,157,455.10 & 1,143,326.29 \\  \hline
40 & 1 & 5 & 7.05 & 0.00 & 4,183.40 & 1,471,754.78 & 1,448,886.60 \\ 
40 & 2 & 5 & 46.27 & 0.00 & 38,245.80 & 1,489,207.54 & 1,470,915.90 \\ 
40 & 3 & 4 & 2,169.68 & 0.04 & 578,702.00 & 1,468,933.86 & 1,450,591.86 \\  \hline
50 & 1 & 5 & 109.11 & 0.00 & 45,634.00 & 1,616,374.96 & 1,578,331.35 \\ 
50 & 2 & 5 & 647.97 & 0.00 & 356,144.00 & 1,738,233.62 & 1,699,686.21 \\ 
50 & 3 & 4 & 2,227.34 & 0.14 & 640,487.00 & 1,659,648.07 & 1,629,222.31 \\  \hline
60 & 1 & 3 & 5,685.70 & 0.12 & 1,032,625.60 & 2,046,328.13 & 1,998,597.77 \\ 
60 & 2 & 4 & 2,164.00 & 0.04 & 370,454.00 & 2,210,560.61 & 2,167,466.49 \\ 
60 & 3 & 4 & 2,373.70 & 0.09 & 512,938.60 & 1,969,974.84 & 1,934,643.82 \\     \bottomrule
    \end{tabular}
\end{table}

In Table~\ref{tab:3index_oneDepot} we display results 
for a single depot, solved with (\texttt{VI}) valid inequalities. In particular, we present the number (\#opt) of instances that were solved to optimality (out of 10), the average CPU times ($t$) in seconds, the average percentage gap at the time limit (3 hours) of all instances, the average number of B\&B nodes (\#nB\&B), the average best found objective function value ($z^*$) at termination and the average lower bound at the root node ($LB_{R}$).

The results in Table~\ref{tab:3index_oneDepot} show that instances with a single-depot EVSP become more difficult with an increasing number of service trips, but not necessarily for an increasing number of charging stations, as the running times for the same number of service trips are sometimes higher with fewer charging stations. Furthermore, Table~\ref{tab:3index_oneDepot} shows that our 
formulation is able to solve nearly all single-depot instances with up to 40 service trips to proven optimality, and most of the instances with 50 and 60 service trips, demonstrating the effectiveness of the formulation for a single depot.

\subsection{Results for a low number of depots}
\label{sec:comp:low}
Next, we consider a set of 270 instances with two, three, and four depots, so $|K| = |O| = |D| \in \{2,3,4\}$.
For each value of $|K|$ we investigate five instances for each combination of $|V^I| \in \{10, 20, 30, 40, 50, 60\}$ service trips and $|C| \in \{1,2,3\}$ charging stations. 

\begin{table}[!htb]
    \centering
    \caption{Results for all formulations for a low number of depots.}
    \begin{scriptsize}
    
    \label{tab:twoThreeDepots}
    \begin{tabular}{lll|rrrrrrrr}
    \toprule
        setting & &  & \#opt & $t$(s) & gap(\%) & \#nB\&B & $z^*$ & $LB_R$ & \#cuts & $t_{cut}$(s) \\ \midrule
        
        \texttt{3i+VI} & & & 175 & 4,287.84 & 0.22 & 400,084.60 & 1,265,436.67 & 1,232,878.84 & ~ &  \\ \hline
        \texttt{2i-IP+VI} & \texttt{I} & \texttt{One} & 196 & 3,353.96 & 0.14 & 529,983.00 & 1,265,443.29 & 1,232,935.81 & 498.80 & 5.08 \\ 
        ~ & \texttt{I} & \texttt{All} & 206 & 2,966.02 & 0.12 & 443,823.40 & 1,265,347.58 & 1,232,880.87 & 299.49 & 1.77 \\ 
        ~ & \texttt{IF} & \texttt{One} & 199 & 3,183.39 & 0.15 & 163,507.80 & 1,265,500.31 & 1,232,976.81 & 122,771.76 & 1,001.90 \\ 
        ~ & \texttt{IF} & \texttt{All} & 196 & 3,392.24 & 0.53 & 162,428.40 & 1,265,598.42 & 1,232,796.79 & 364,862.30 & 1,262.44 \\ \hline
        \texttt{2i-CC+VI} & \texttt{I} & \texttt{One} & 202 & 3,089.69 & 0.14 & 412,167.30 & 1,265,398.67 & 1,232,943.92 & 323.79 & 6.16 \\ 
        ~ & \texttt{I} & \texttt{All} & 206 & 2,995.78 & 0.14 & 385,027.70 & 1,265,393.44 & 1,232,887.54 & 226.26 & 4.07 \\ 
        ~ & \texttt{IF} & \texttt{One} & 141 & 5,423.97 & 4.84 & 3,974.10 & 1,248,703.81 & 1,233,114.30 & 30,155.29 & 302.43 \\ 
        ~ & \texttt{IF} & \texttt{All} & 141 & 5,335.53 & 4.88 & 4,677.80 & 1,245,562.88 & 1,233,046.19 & 25,999.09 & 189.79 \\  \bottomrule
    \end{tabular}
    \end{scriptsize}
\end{table}

In Table~\ref{tab:twoThreeDepots} we display average results over all 270 instances for each available setting, namely for using 
\begin{itemize}
    \item the 3-index formulation~\eqref{z3_MDEVSP_model} (named \texttt{3i}) with (\texttt{VI}) valid inequalities ~\eqref{val_i:z3_MDEVSP_lowerboundbuses},~\eqref{val_i:z3_MDEVSP_consumeEnergyST},~\eqref{val_i:z3_MDEVSP_loadEnergy}, and~\eqref{val_i:z3_MDEVSP_energyDestination}, 
    \item the 2-index formulation~\eqref{z2_MDEVSP_model}, which includes the infeasible path constraints~\eqref{constr:z2_MDEVSP_infeasiblePath} (named \texttt{2i-IP}), with separating only integer (\texttt{I}) or both integer and fractional (\texttt{IF}) solutions, with including all (\texttt{All}) or only the first (\texttt{One}) infeasible path constraint, with (\texttt{VI})  valid inequalities~\eqref{val_i:z2_MDEVSP_lowerboundbuses},~\eqref{val_i:z2_MDEVSP_consumeEnergyST},~\eqref{val_i:z2_MDEVSP_loadEnergy}, and~\eqref{val_i:z2_MDEVSP_energyDestination}, and
    \item the 2-index formulation~\eqref{z2_MDEVSP_model} without~\eqref{constr:z2_MDEVSP_infeasiblePath} but with the connectivity constraints~\eqref{val_i:z2_MDEVSP_maxFlow} (named \texttt{2i-CC}), with separating only integer (\texttt{I}) or both integer and fractional (\texttt{IF}) solutions, with including all (\texttt{All}) or only the first (\texttt{One}) connectivity constraint, and
    with (\texttt{VI}) valid inequalities~\eqref{val_i:z2_MDEVSP_lowerboundbuses},~\eqref{val_i:z2_MDEVSP_consumeEnergyST},~\eqref{val_i:z2_MDEVSP_loadEnergy}, and~\eqref{val_i:z2_MDEVSP_energyDestination}. 
\end{itemize}
For each of these settings, we give the number (\#opt) of instances that were solved to optimality (out of 270), the average CPU times ($t$) in seconds, the average percentage gap at the time limit (3 hours) of all instances, the average number of B\&B nodes (\#nB\&B), the average best found objective function value ($z^*$) at termination, the average lower bound at the root node ($LB_{R}$), the average number of cuts (\#cuts, i.e., added infeasible path constraints~\eqref{constr:z2_MDEVSP_infeasiblePath} or added connectivity constraints~\eqref{val_i:z2_MDEVSP_maxFlow}) and the average CPU time ($t_{cut}$) in seconds for the separation. For instances where no feasible solution was found within the time limit, we assume a gap of 100\% and exclude $z^*$ from the calculations.

The results in Table~\ref{tab:twoThreeDepots} clearly show that for \texttt{2i-CC+VI} both versions separating fractional and integer solutions \texttt{IF} perform much worse and solve much fewer instances to optimality within the time limit, demonstrating the superiority of \texttt{I} for \texttt{2i-CC+VI}.

For \texttt{2i-IP+VI} the picture is not that clear, as in the setting \texttt{One}, \texttt{IF} performs slightly better than \texttt{I} with three more solved instances (199 vs.\ 196), while in the setting \texttt{All}, \texttt{IF} solves less instances than \texttt{I} (196 vs. 206).

Furthermore, the 2-index formulation with infeasible path constraints (\texttt{2i-IP+VI}) and both \texttt{I} and \texttt{IF} separations, as well as the 2-index formulation with connectivity cuts (\texttt{2i-CC+VI}) and \texttt{I} separations, outperform the 3-index model in all settings (at least 196 instances are solved to optimality in all \texttt{2i-I+VI} settings, while 175 instances are solved to optimality in \texttt{3i+VI}).
Overall, the best settings of this computational evaluation are \texttt{2i-IP+VI+I+All} and \texttt{2i-CC+VI+I+All} with 206 instances solved to optimality each.

\subsection{Results for a high number of depots}
\label{sec:comp:high}
Finally, we consider a set of instances with $|K| = |O| = |D| \in \{6,8\}$ depots. Here we consider for each value of $|K|$ five instances with $|C|=\{3,4,6\}$ charging stations for each number of service trips $|V^I| \in \{30,40,50,60,70,80\}$, resulting in 180 instances in total.

\begin{table}[!htb]
    \centering
    \caption{Results for 2-index formulations and a high number of depots.}
    \begin{scriptsize}
    \label{tab:2dxHighDepots}
    \begin{tabular}{ll|rrrrrrrrrr}
\toprule
    &  & \#opt & $t$(s) & gap(\%) & \#nB\&B & $z^*$ & $LB_R$ & \#cuts & $t_{cut}$(s) \\ 
    \midrule
        \texttt{2i-IP+VI+I} & \texttt{One} & 46 & 8,420.18 & 38.02 & 536,048.60 & 1,579,065.89 & 1,785,904.08 & 6,224.49 & 92.52 \\ 
        \texttt{2i-IP+VI+I} & \texttt{All} & 60 & 7,478.90 & 14.44 & 408,568.40 & 1,779,073.37 & 1,785,904.24 & 2,903.48 & 7.20 \\ \hline
        \texttt{2i-CC+VI+I} & \texttt{One} & 58 & 7,601.72 & 14.95 & 339,778.70 & 1,768,086.18 & 1,785,903.89 & 2,352.59 & 72.59 \\ 
        \texttt{2i-CC+VI+I} & \texttt{All} & 65 & 7,308.79 & 4.92 & 300,350.80 & 1,832,566.14 & 1,785,905.88 & 1,474.92 & 24.26 \\ 
    \bottomrule
\end{tabular}
\end{scriptsize}
\end{table}

Table~\ref{tab:2dxHighDepots} gives the average results obtained with the 2-index model~\eqref{z2_MDEVSP_model} in the two versions \texttt{2i-IF} and \texttt{2i-CC}, adding either all cuts (\texttt{All}) or only one (\texttt{One}) in each call of the separation routine. They are only separated on integer solutions (setting \texttt{I} from above) and the valid inequalities are active (setting \texttt{VI} from above). The columns of Table~\ref{tab:2dxHighDepots} are defined analogously to the columns of Table~\ref{tab:3index_oneDepot} and Table~\ref{tab:twoThreeDepots}.

The settings \texttt{2i-IP+VI+I} and \texttt{2i-CC+VI+I} in combination with adding \texttt{All} cuts solve 60 and 65 instances out of 180 to optimality, respectively, and clearly outperform both settings with \texttt{One}. When comparing computation times, \texttt{2i-IP+VI+I+All} requires 7,478.90 seconds on average and solves 60 instances to optimality, while \texttt{2i-CC+VI+I+All} requires 7,308.79 seconds on average and solves 65 instances to optimality, indicating that \texttt{2i-CC+VI+I+All} has a clear edge over \texttt{2i-IP} for instances with a high number of depots. Overall, not only the number of optimally solved instances, but also the optimality gap is much better for both settings with \texttt{All} cuts added.
In the end, the setting \texttt{2i-CC+VI+I+All} can be determined as the clear winner in these runs.

\begin{table}[!htb]
    \centering
    \small
    \caption{Detailed results of the 2-index and 3-index formulations for a high number of depots.}
    \label{tab:allHighDepots}
 \begin{tabular}{ll|rrr|rrr|rrr}
 \toprule
\multicolumn{2}{l|}{} & \multicolumn{3}{l|}{2-index (\texttt{2i-IP+VI+I+All})} & \multicolumn{3}{l|}{2-index (\texttt{2i-CC+VI+I+All})} & \multicolumn{3}{l}{3-index (\texttt{3i+VI})}\\
\midrule
$|V^I|$ & $|K|$ & \#opt & $t$(s) & gap(\%) & \#opt & $t$(s) & gap(\%) & \#opt  &  $t$(s) & gap(\%) \\
\midrule
        30 & 6 & 15 & 311.64 & 0.00 & 15 & 142.63 & 0.00 & 9 & 4,617.07 & 0.23 \\ 
        30 & 8 & 15 & 312.54 & 0.00 & 15 & 113.45 & 0.00 & 10 & 4,143.91 & 0.35 \\ \hline
        40 & 6 & 10 & 4,088.86 & 0.12 & 11 & 3,540.53 & 0.07 & 6 & 6,830.38 & 0.43 \\ 
        40 & 8 & 11 & 3,912.17 & 0.08 & 11 & 4,199.32 & 0.08 & 4 & 8,433.61 & 0.88 \\  \hline
        50 & 6 & 6 & 6,654.56 & 0.44 & 6 & 6,599.81 & 0.45 & 4 & 9,272.83 & 1.04 \\ 
        50 & 8 & 3 & 9,573.76 & 0.54 & 6 & 8,441.46 & 0.43 & 1 & 10,707.20 & 1.81 \\  \hline
        60 & 6 & 0 & 10,803.65 & 7.34 & 0 & 10,803.51 & 0.84 & 0 & 10,803.86 & 1.36 \\ 
        60 & 8 & 0 & 10,803.75 & 1.07 & 1 & 10,622.99 & 0.93 & 0 & 10,803.65 & 1.95 \\ 
\bottomrule
\end{tabular}
\end{table}

In Table~\ref{tab:allHighDepots}, we compare the best performing 2-index-based B\&C algorithm \texttt{2i-CC+VI+I+All} for instances with a high number of depots with the 2-index approach \texttt{2i-IP+VI+I+All} and with
solving the 3-index formulation with valid inequalities \texttt{3i+VI} with CPLEX directly for
120 out of the 180 instances in detail, differentiating between the number of service trips $|V^I|$ and the number of depots $|K|$. For both the 2- and the 3-index formulations, we present the number (\#opt) of instances that were solved to optimality (out of 15), the average percentage gap at the time limit (3 hours) of all instances, and the average CPU times ($t$) in seconds. We show
that the 2-index formulation \texttt{2i-CC+VI+I+All} demonstrates superior overall performance in terms of number of instances solved to optimality (65 vs. 60 vs. 34). Furthermore, the 2-index formulation consistently achieves significantly lower average optimality gaps, and it generally requires notably less computational time, demonstrating the efficiency of our 2-index approach with connectivity constraints for a high number of depots.

\begin{table}[!htb]
    \centering
    \small
    \caption{Comparison of lower bounds for the 2-index and 3-index formulation, applied to instances with a high number of depots and 70 or 80 service trips.}
    \label{tab:STHighDepots}
 \begin{tabular}{cc|cc|cc}
 \toprule
\multicolumn{2}{c|}{} & \multicolumn{2}{c|}{2-index (\texttt{2i-CC+VI+I+All})} & \multicolumn{2}{c}{3-index (\texttt{3i+VI})}\\
\midrule
$|V^I|$ & $|K|$ & $LB_R$ & $LB$ & $LB_R$ & $LB$ \\
\midrule
        70 & 6 & 2,181,500.90 & 2,230,294.71 & 2,181,574.43 & 2,204,635.79 \\ 
        70 & 8 & 2,168,404.87 & 2,220,856.65 & 2,168,201.98 & 2,178,579.94 \\ \hline
        80 & 6 & 2,494,965.11 & 2,548,422.41 & 2,495,026.01 & 2,511,729.53 \\ 
        80 & 8 & 2,381,890.18 & 2,433,957.84 & 2,381,694.60 & 2,394,358.14 \\  
\bottomrule
\end{tabular}
\end{table}

Table~\ref{tab:STHighDepots} compares the performance of the 2-index (\texttt{2i-CC+VI+I+All}) and 3-index (\texttt{3i+VI}) formulations on the remaining 60 instances with a high number of depots and 70 or 80 service trips. These specific instances are shown because neither formulation was able to solve them to optimality within the time limit. As a result, we focus on comparing the quality of the lower bounds obtained by each model. For each case, we report the lower bound at the root node ($LB_R$) and the best lower bound obtained within the time limit ($LB$). The results show that both formulations yield similar lower bounds at the root node. However, the 2-index formulation generally achieves slightly stronger bounds at the time limit compared to the 3-index model, indicating a more effective exploration of the solution space.

\section{Comparison of propulsion technologies}
\label{techComparison}

In this section, we
compare different propulsion technologies, namely diesel (DB), battery-electric (BEB), and fuel cell-electric buses (FCEB), by evaluating their impact on vehicle scheduling.
A key aspect of this analysis is to determine the number of buses required for each propulsion type, which is crucial for planning fleet transitions. In the following, first, we describe the computational setup and the attributes of our newly generated test instances in Section~\ref{rl:sec:comp:instances}. Next, the obtained results are discussed in Section~\ref{rl:sec:comp:base}. We examine results under cold temperatures in Section~\ref{rl:sec:comp:cold} and for a scenario preserving battery life in Section~\ref{rl:sec:comp:battery}.

\subsection{Computational setup and test instances}
\label{rl:sec:comp:instances}

In our computational study we use the same computational setup as described in Section~\ref{sec:comp:instances}.
We further generate a second set of realistic instances, designed to reflect real-world behavior based on data provided by our project partners. This data includes observed patterns of service trips as well as their spatial and temporal distribution.
In addition, these realistic instances incorporate vehicle and refueling/recharging specifications for three distinct types of buses, them being DBs, BEBs, and FCEBs, reflecting a range of technological and operational characteristics. 
In total, we consider 34 test instances that vary in size and structure to reflect a broad range of operational settings. These include 10 base and 24 combined instances, where base instances consist of one bus line each. The combined instances are created by randomly merging base instances in different patterns to reflect various ways of grouping the bus lines. The time limit for solving each instance is set to three hours.
In the following, the generation of base instances is described in Section~\ref{rl:sec:base}, while the creation of combined instances is detailed in Section~\ref{rl:sec:combined}.

\subsubsection{Base instances}
\label{rl:sec:base}

Each of the 10 base instances accounts for one separate bus line and is defined by $|K| = |O| = |D| = |C| = 1$ depot and charging station. For the number of service trips, three instances are generated for each value in $|V^I| \in \{10, 20, 30\}$, and one instance for $|V^I| = 40$.
The BEB charging stations are co-located with the depot sites. 
Because high capital and operating costs, strict safety regulations, and the complexity of hydrogen supply logistics make it impossible to install a hydrogen refueling station at every depot location, FCEBs must refuel at a single off-site hydrogen station, the location of which is the same in all instances. It is assigned randomly within the 50 km × 50 km square in the Euclidean plane, following a uniform distribution. The square size is chosen similar to the size of the project's model regions.
For DBs, refueling is not allowed during the day as it is not necessary for our project partners. Therefore, there is no refueling stations for DBs.

Each service trip going in one direction has a corresponding return trip in the opposite direction. We refer to these service trips as forward and backward trips, respectively. 
In each instance, a line has a generally fixed start and end location $\ell_i^s$ and $\ell_i^e$
for forward trips, with occasional variations, which occur in 3\% to 17\% of the service trips as observed in the project data. 
The default start location of a line $\ell_i^s$, as well as the locations of depots $\ell_k$ for each $k \in K$, are again randomly assigned within a 50 km × 50 km square in the Euclidean plane, following a uniform distribution. This random assignment determines the latitude and longitude coordinates for each location. Note that for the whole instance set, we only define three distinct depot locations.

The end location of each forward trip $\ell_i^e$ is placed at a minimum distance of 15 km from its start location. Insights from project data indicate that most lines are approximately 15 km long, with fewer instances as the distance increases, though some outliers extend up to 60 km. To reflect this distribution, we sample distances from a normal distribution with a mean of 15 km and a standard deviation of 15 km, and then reject values below 15 km and resample, until a distance greater or equal to 15 is obtained. Once a target distance is determined, we generate 200 candidate points in the square and select the one closest to the target distance from the start location as the end location of the line.
For modified service trips, the start and/or end locations are adjusted. Let $d_{i}$ denote the Euclidean distance between the default start and end locations $\ell_i^s$ and $\ell_i^e$ of a line. 
The new start and/or end locations are randomly selected within a radius of $d_i/2$ km from their original positions.
From these locations we compute the distances $d_{ij}$ as the Euclidean distances between the (end) location of $i$ and the (start) location of $j$ for all $i,j \in V^I \cup O \cup D \cup C$. Moreover, we compute the energy usage $p_{ij}$ as $p_{ij} = \theta d_{ij}$.
For the travel times $t_{ij}$ (in minutes) for deadhead trips, we assume an average vehicle speed of 60 km per hour, equivalent to 1 km per minute, so $t_{ij} = d_{ij}$.

The duration $u_i$ of a line is randomly determined to be between $1.7$ and $3$ times the distance $d_i$ of the line. 
For each service trip $i \in V^I$, the duration is then randomly generated within the range of $0.9u_i$ to $1.1u_i$. In the project instances, most service trips have the same start and end stops, but few of them (3 to 17 \%) are shorter or longer, with different start and/or end points.
To imitate this behavior, we modify between 3 and 17 \% of the service trips.
For modified service trips, the distance is recalculated with the new start and/or end points, and the duration $u_i$ is adjusted accordingly using the same coefficient. Then, the durations of forward and backward trips are generated independently between $0.9u_i$ to $1.1u_i$, i.e., they may have different durations. The energy usage $q_i$ of each service trip $i \in V^I$ is calculated as $q_i = \theta u_i$.

Like in \cite{Carpaneto1989}, the start time $s_i$ (in minutes since midnight) for each forward trip $i \in V^I$ is generated as a random integer in the intervals [420, 480] and [1020, 1080] with a probability of 15\% each. In our case, to increase the likelihood of service trips occurring before 7 a.m. and after 6 p.m., we also generate $s_i$ in the intervals [300, 420], [480, 1020], and [1080, 1410] with a combined probability of 70\%.
Moreover, we compute the end time $e_i$ of each service trip $i \in V^I$ as $e_i = s_i + u_i$. 

The start time of the corresponding backward trip $s_j$ for $j \in V^I$ is given by the end time of the forward trip $e_i$ plus an additional 10 minutes of idle time at the end location of the forward trip $i \in V^I$. Specifically, 
$s_j = e_i + 10$, where 10 minutes represents idle time. 
For each depot index $k \in K$ we generate the number of vehicles available $b_k$ in the same way as for our benchmark instances, namely as a uniformly random integer in $[3 + \frac{1}{3|K|} V^I,
3 + \frac{1}{2|K|}|V^I|]$.

The parameters related to vehicles for DBs are set to $\theta= 0.47$ \citep{chikishev2022assessment}, 
and $s^{max}=350$ \citep{li2019mixed}.

Further, we define the parameters for FCEBs and their fueling technology. We set $s^{max}=37$ and $\theta= 0.08$ \citep{reithuber2025energy}.
According to our project partners, fully refueling a fuel cell bus takes 12 minutes, therefore we set $t^{max} = t^{min} = 12$. Note, that this results in $r=s^{max}/t^{max} \approx 3.08$.
Since the refueling station for FCEBs is not located at the depot, we set the minimum state of charge (SOC) at any node except the depot to $s^{min}=0.2s^{max} = 7.4$. Furthermore, we define the minimum required SOC when returning to the depot as $s^{min}_{dep}=\lfloor s^{max}-  \max\limits_{a \in C,\ d \in V^D} p_{ad}\rfloor $. 

For BEBs, we use specifications based on the MAN Lion’s City 12 E. 
According to \cite{MAN_LionsCity12E_datasheet}, the range of BEBs under favorable operating conditions, described as low operational demands in terms of average speed, loading, topography, and air-conditioning, is up to 380 km. Since we assume more demanding conditions, we set the range of our BEBs to 80\% of the given one, resulting in 304 km. We set $s^{max}=445$ and further use said range to calculate $\theta= s^{max}/304 = 1.46$, where we divide $s^{max}$ by the range in km. Based on feedback from our project partners, we assume a charger power of 130 kW, therefore $r=130/60 \approx 2.17$.
Note, that this yields $t^{max} = (s^{max}-s^{min})/r=164$  minutes. 
Furthermore, since we assume depots are equipped with recharging stations, we fix $s^{min}_{dep}=s^{min}=0.2s^{max} = 89$ and $t^{min}=15$. An overview of the attributes and other parameters specified for each bus type for our standard setting is given in Table~\ref{tab:base:busparam}.

\begin{table}[H]
\centering
\caption{Parameter specification for FCEBs and BEBs (standard setting).}
\label{tab:base:busparam}
\begin{tabular}{lll}
\toprule
& FCEB  & BEB \\
\midrule
$s^{max}$  & $37$ kg\textsuperscript{a} & $445$ kWh\textsuperscript{b} \\
$s^{min}$ & 7.4 kg & 89 kWh \\
$s^{min}_{dep}$ & $ \lfloor s^{max}- \max\limits_{a \in C,\ d \in V^D} p_{ad} \rfloor$ & 89 kWh \\
$\theta$ & 0.08 kg/km\textsuperscript{a} & 1.46 kWh/km\textsuperscript{b} \\
$t^{max}$  & 12 min\textsuperscript{c} & 164 min \\
$t^{min}$  & 12 min & 15 min\\
$r$  & 3.08 kg/min & $2.17$ kW/min\textsuperscript{c} \\
\bottomrule
\end{tabular}
\end{table}

\vspace{0.5em}
\scriptsize
\noindent
\textsuperscript{a} From \cite{reithuber2025energy}.\\
\textsuperscript{b} From \cite{MAN_LionsCity12E_datasheet}.\\
\textsuperscript{c} From project partner and field observations.
\normalsize

\subsubsection{Combined instances}
\label{rl:sec:combined}

The 24 combined instances are generated by merging base instances to increase network size and complexity while preserving the structure of each individual problem. Specifically, we begin by randomly merging two base instances at a time until we obtain ten combined instances with two lines. These are then extended by one randomly selected line to create ten instances with three lines.  Finally, the four instances with four lines are constructed by merging two of the two-line instances and by extending two of the three-line instances with an additional line.
By grouping together base instances with one line each, we ensure that each combined instance retains the same underlying topology and parameters, but varies systematically in size. This systematic aggregation allows us to assess how solution methods scale when the number of lines or depots grows, and to observe interactions between different lines that would not appear if each base instance were solved separately. It also produces instances whose data volumes more closely resemble real‐world systems, without introducing fundamentally new network patterns.

When multiple lines are merged, shared depot and refueling resources can be utilized more efficiently, smoothing peak loads and reducing idle time. In separate planning (i.e., the base instances), each line is assigned enough buses to independently cover its own peak load and refueling needs, which can result in underutilized vehicles during off-peak periods.  In the combined instances, buses can be flexibly assigned across lines, and because recharging is permitted at every depot, deadhead trips to and from depots are minimized. This coordinated scheduling approach can lead to a reduction in the total fleet size while providing valuable managerial insights into how resources can be shared more efficiently across different bus lines.

Thus, each combined instance has a number of service trips $|V^I| \in \{30, 40, \dots, 100\}$ and
$|K| = |O| = |D| = |C| \in \{1,2, 3\}$ depots and charging stations. From the base instances, we inherit the locations of depots and, therefore, charging stations as well as service trips. It follows that only one charging station is deployed for FCEBs.

\subsection{Computational results}
\label{rl:sec:comp:base}

We first solve each of the base instances individually, treating each as a distinct operational line. Subsequently, we solve the combined instances, which are constructed by merging multiple base instances. These multi-line cases include a broader set of service trips and greater operational diversity, allowing us to analyze the benefits of integrated planning across lines while maintaining consistency in the underlying problem structure.

\begin{table}[ht]
\centering
\begin{scriptsize}
\caption{Results for the base instances for each technology.}
\begin{tabular}{cc|crc|ccrc|ccrc}
\toprule
\multicolumn{2}{c|}{}
& \multicolumn{3}{c|}{DB} 
& \multicolumn{4}{c|}{BEB} 
& \multicolumn{4}{c}{FCEB} \\
line & $|V^I|$ & \#b & \multicolumn{1}{c}{$z^*$} & gap(\%) 
  & \#b & \#c & \multicolumn{1}{c}{$z^*$} & gap(\%) 
  & \#b & \#c & \multicolumn{1}{c}{$z^*$} & gap(\%) \\
\midrule
1 & 10 & 2 &  200,048.28  & 0.00 & 3 & 2 &  308,316.92  & 0.00 & 2 & 3 &  212,020.22  & 0.00 \\ 
2 & 10 & 2 &  200,074.08  & 0.00 & 2 & 1 &  204,341.18  & 0.00 & 2 & 2 &  208,013.24  & 0.00 \\ 
3 & 10 & 2 &  200,039.52  & 0.00 & 2 & 2 &  208,287.46  & 0.00 & 2 & 2 &  208,011.38  & 0.00 \\ \hline
4 & 20 & 3 &  300,145.02  & 0.00 & 4 & 2 &  408,850.98  & 0.00 & 3 & 3 &  312,024.69  & 0.00 \\ 
5 & 20 & 3 &  300,055.02  & 0.00 & 4 & 4 &  416,468.49  & 0.00 & 3 & 5 &  320,010.77  & 0.00 \\ 
6 & 20 & 4 &  400,113.50  & 0.00 & 4 & 1 &  404,425.88  & 0.00 & 4 & 4 &  416,023.87  & 0.00 \\ \hline
7 & 30 & 4 &  400,083.90  & 0.00 & 5 & 5 &  520,777.74  & 1.14 & 4 & 6 &  424,023.23  & 1.89 \\ 
8 & 30 & 4 &  400,113.12  & 0.00 & 4 & 2 &  408,435.72  & 1.98 & 4 & 4 &  416,025.28  & 2.88 \\ 
9 & 30 & 5 &  500,162.16  & 0.00 & 6 & 11 &  644,937.16  & 0.00 & 5 & 7 &  528,028.86  & 0.00 \\ \hline
10 & 40 & 6 &  600,272.69  & 0.00 & 10 & 16 &  1,065,785.48  & 0.00 & 6 & 14 &  656,064.36  & 0.00 \\ 
\bottomrule
\end{tabular}
\label{tab:tech_base_base_results}
\end{scriptsize}
\end{table}

Table \ref{tab:tech_base_base_results} presents the computational results for the base instances across all considered propulsion technologies. Displayed are the total number of buses (\#b), total number of charging events (\#c), the percentage gap at the time limit of three hours, and the best found objective function value ($z^*$).

The used setting 
(\texttt{2i+All+VI}) successfully solves to proven optimality all instances for DBs, for BEBs and FCEBs two of the instances with 30 service trips remain with an optimality gap. DB and FCEBs use the same amount of buses for each instance, however, each FCEBs needs to use the charging station at least once before returning to the depot at the end of its schedule. 35 DB or FCEBs, respectively, are needed, to cover all service trips. Meanwhile, the best found schedule for BEBs calls for 44 buses. In their schedule, 46 charging events happen. The big difference in the obligatory buses to cover the schedule between BEBs and the other two technologies diesel and fuel cell-electric can be explained by the fact that the energy requirements of the individual service trips are quite high. In addition to that, we reduced $s^{max}$ for the realistic instances to a more realistic number.
The number of recharges tends to increase with the instance size.
As expected, the required number of buses generally increases as well with the number of trips. The difference of buses needed within each instance size may be due to the placement of the depot location, i.e., whether it is closer or further away from the service trips.

\begin{table}[ht]
\centering
\scriptsize
\caption{Results for all instances solved separately vs. combined, grouped by technology.}
\begin{tabular}{c r l| 
                *{2}{c}| 
                *{4}{c}| 
                *{4}{c}}
\toprule
& & 
  & \multicolumn{2}{c|}{DB} 
  & \multicolumn{4}{c|}{BEB} 
  & \multicolumn{4}{c}{FCEB} \\
  &  &  base
  & \multicolumn{1}{c}{ separate} & \multicolumn{1}{c|}{ combined}
  & \multicolumn{2}{c}{ separate} & \multicolumn{2}{c|}{ combined}
  & \multicolumn{2}{c}{ separate} & \multicolumn{2}{c}{ combined} \\
$|K|$ & $|V^I|$ & lines
   & \#b &  \#b 
   & \#b &  \#c &  \#b &  \#c 
  &  \#b &  \#c &  \#b &  \#c \\
\midrule
1 & 50 & 1, 10 & 8 & 8 & 13 & 18 & 12 & 18 & 8 & 17 & \textit{8} & \textit{17}\\
1 & 30 & 1, 4 & 5 & 5 & 7 & 4 & 7 & 3 & 5 & 6 & 5 & 6\\
2 & 40 & 2, 7 & 6 & 5 & 7 & 6 & \textit{5} & \textit{8} & 6 & 8 & \textit{5} & \textit{7}\\
2 & 40 & 2, 9 & 7 & 7 & 8 & 12 & 7 & 10 & 7 & 9 & \textit{7} & \textit{8}\\
1 & 40 & 3, 9 & 7 & 7 & 8 & 13 & 8 & 10 & 7 & 9 & \textit{7} & \textit{9}\\
1 & 50 & 4, 7 & 7 & 6 & 9 & 7 & \textit{7} & \textit{10} & 7 & 9 & \textit{6} & \textit{8}\\
2 & 60 & 5, 10 & 9 & 8 & 14 & 20 & \textit{12} & \textit{20} & 9 & 19 & 8 & 18\\
1 & 50 & 6, 9 & 9 & 7 & 10 & 12 & \textit{8} & \textit{10} & 9 & 11 & \textit{7} & \textit{10}\\
2 & 60 & 7, 8 & 8 & 8 & 9 & 7 & \textit{8} & \textit{8} & 8 & 10 & \textit{8} & \textit{9}\\
2 & 60 & 8, 9 & 9 & 6 & 10 & 13 & \textit{7} & \textit{11} & 9 & 11 & \textit{6} & \textit{10}\\
 \midrule 
 \multicolumn{3}{r|}{\text{\diameter}} & 7.5 & 6.7 & 9.5 & 11.2 & 8.1 & 10.8 & 7.5 & 10.9 & 6.7 & 10.2 \\ 
 \midrule 
2 & 70 & 1, 10, 6 & 12 & 9 & 17 & 19 & \textit{11} & \textit{21} & 12 & 21 & \textit{9} & \textit{20}\\
1 & 70 & 1, 4, 10 & 11 & 11 & 17 & 20 & 15 & 20 & 11 & 20 & \textit{11} & \textit{19}\\
3 & 70 & 2, 7, 9 & 11 & 10 & 13 & 17 & \textit{14} & \textit{15} & 11 & 15 & \textit{10} & \textit{15}\\
2 & 70 & 2, 9, 8 & 11 & 8 & 12 & 14 & \textit{10} & \textit{11} & 11 & 13 & \textit{8} & \textit{13}\\
2 & 60 & 3, 9, 5 & 10 & 10 & 12 & 17 & \textit{10} & \textit{14} & 10 & 14 & \textit{10} & \textit{13}\\
2 & 80 & 4, 7, 8 & 11 & 9 & 13 & 9 & \textit{10} & \textit{11} & 11 & 13 & \textit{9} & \textit{13}\\
2 & 70 & 5, 10, 2 & 11 & 10 & 16 & 21 & \textit{14} & \textit{25} & 11 & 21 & \textit{10} & \textit{19}\\
2 & 80 & 6, 9, 7 & 13 & 10 & 15 & 17 & \textit{12} & \textit{16} & 13 & 17 & \textit{10} & \textit{15}\\
2 & 70 & 7, 8, 1 & 10 & 8 & 12 & 9 & \textit{9} & \textit{12} & 10 & 13 & \textit{8} & \textit{13}\\
3 & 90 & 8, 9, 7 & 13 & 10 & 15 & 18 & - & - & 13 & 17 & - & -\\
 \midrule 
 \multicolumn{3}{r|}{\text{\diameter}} & 11.3 & 9.5 & 14.2 & 16.1 & 11.7 & 16.1 & 11.3 & 16.4 & 9.4 & 15.6 \\ 
 \midrule 
3 & 100 & 1, 10, 6, 8 & 16 & 10 & 21 & 21 & - & - & 16 & 25 & - & -\\
2 & 80 & 1, 4, 6, 9 & 14 & 11 & 17 & 16 & \textit{19} & \textit{18} & 14 & 17 & - & -\\
3 & 80 & 2, 7, 3, 9 & 13 & 12 & 15 & 19 & - & - & 13 & 17 & - & -\\
3 & 90 & 4, 7, 8, 3 & 13 & 9 & 15 & 11 & \textit{17} & \textit{26} & 13 & 15 & - & -\\
 \midrule 
 \multicolumn{3}{r|}{\text{\diameter}} & 14.0 & 10.5 & 17.0 & 16.8 & 18.0 & 22.0 & 14.0 & 18.5 & - & - \\ 
\bottomrule
\end{tabular}
\label{tab:base:2dxintegrated_separate}
\end{table}

Table \ref{tab:base:2dxintegrated_separate} presents the total average number of buses (\#b) and charging events (\#c) across DB, BEB, and FCEB technologies, under two planning approaches. In the separate case, each base line is solved individually, and the total number of buses and charging events is obtained by summing up the results across lines. In the combined case, the same set of lines is solved together as a combined instance, which allows for shared use of buses and potentially more efficient charging schedules. All instances not solved to optimality are displayed in italic. For instances, where we were not able to find a feasible solution, ``-'' is displayed.
We use the setting (\texttt{2i-IP+All+I+VI}) for solving all instances with more than one depot from now on.
Each row corresponds to a specific combined instance, and we present the the number of depots $|K|$, the number of service trips $|V^I|$ and which base lines are combined.

The results for two lines show that, as expected, solving the lines jointly tends to reduce the number of buses slightly for DBs, BEBs, and FCEBs.
Looking at all technologies, in almost all cases, following the combined approach spares at least one bus compared to serving the two lines individually, even if some of the combined instances are not solved to optimality.
The computed averages confirm this trend. Although fewer buses are required, the total number of charging events declines only slightly, implying that each bus must recharge more frequently on average.

We extend the analysis to scenarios involving three and four lines. As the number of lines increases, so does the complexity of the problem, which is reflected in the higher number of required buses and charging events, especially for BEBs. We observe that combining three lines consistently reduces the number of buses by around two buses across all technologies, even though many instances are not solved to proven optimality.
The largest fleet reduction for BEBs is observed when solving lines 1, 10, 6 combined instead of separately, yielding a saving of six buses. In comparison, the same lines can be served with only nine FCEBs or DBs in total when solved as a combined instance. Regarding charging events, the results for BEBs are mixed: depending on the instance, the number of events may decrease, remain unchanged, or even increase. For FCEBs, by contrast, the number of charging events tends to slightly decrease or remain unchanged. Since less buses are serving the lines, this implies a higher charging frequency of the remaining buses.
Neither for BEB, nor for FCEB, we could find a feasible solution in the combined case for line 8, 9, 7, reflecting the increasing computational difficulty as more lines are integrated and the number of service trips rises.

For DBs, the benefits of combining four lines remain evident: the average fleet size decreases from 14.0 to 10.5 buses. In contrast, no feasible solutions could be obtained for two BEB instances and for all FCEB instances, indicating that solving larger integrated problems for electric fleets is considerably more demanding computationally. 
For BEBs with four combined lines, the reported fleet sizes are in some cases even increasing, which is explained by very large optimality gaps: for instance, the case with lines 1, 4, 6, 9 has a gap of 44.25\%, while the case with lines 4, 7, 8, 3 shows a gap of 50.16\%, indicating that the solutions are likely far from the true optimum.

Across all instances, DBs and FCEBs require the same number of buses, regardless of whether the lines are planned separately or in combination. In general, the results reinforce that integrating multiple lines into a joint optimization problem leads to reductions in fleet size.

\begin{table}[ht]
\centering
\scriptsize
\caption{Average percentage gaps and CPU times ($t$) in seconds across all technologies, solved combined.}
\begin{tabular}{c r l| 
                *{2}{r}| 
                *{2}{r}| 
                *{2}{r}}
\toprule
& & 
  & \multicolumn{2}{c|}{DB} 
  & \multicolumn{2}{c|}{BEB} 
  & \multicolumn{2}{c}{FCEB} \\
  $|K|$ & $|V^I|$ & lines
   & gap(\%) &  $t$(s) 
    & gap(\%) &  $t$(s) 
     & gap(\%) &  $t$(s) \\
\midrule
1 & 50 & 1, 10 & 0.00 & 16.60 & 0.00 & 12.05 & 1.70 & TL\\
1 & 30 & 1, 4 & 0.00 & 8.71 & 0.00 & 130.43 & 0.00 & 65.87\\
2 & 40 & 2, 7 & 0.00 & 8.70 & 1.05 & TL & 1.51 & TL\\
2 & 40 & 2, 9 & 0.00 & 8.74 & 0.00 & 334.05 & 0.55 & TL\\
1 & 40 & 3, 9 & 0.00 & 8.65 & 0.00 & 348.90 & 1.09 & TL\\
1 & 50 & 4, 7 & 0.00 & 8.91 & 1.64 & TL & 2.53 & TL\\
2 & 60 & 5, 10 & 0.00 & 78.62 & 1.55 & TL & 0.00 & 118.55\\
1 & 50 & 6, 9 & 0.00 & 8.87 & 1.41 & TL & 1.62 & TL\\
2 & 60 & 7, 8 & 0.00 & 9.61 & 3.87 & TL & 4.31 & TL\\
2 & 60 & 8, 9 & 0.00 & 9.62 & 18.81 & TL & 6.25 & TL\\
 \midrule 
 \multicolumn{3}{r|}{\text{\diameter}} & 0.00 & 16.70 & 2.83 & 6,568.24 & 1.96 & 8,665.96 \\ 
 \midrule 
2 & 70 & 1, 10, 6 & 0.00 & 13.39 & 13.55 & TL & 8.17 & TL\\
1 & 70 & 1, 4, 10 & 0.00 & 10.36 & 0.00 & 854.35 & 3.40 & TL\\
3 & 70 & 2, 7, 9 & 0.00 & 9.02 & 31.55 & TL & 5.66 & TL\\
2 & 70 & 2, 9, 8 & 0.00 & 9.34 & 23.39 & TL & 6.11 & TL\\
2 & 60 & 3, 9, 5 & 0.00 & 9.08 & 2.91 & TL & 1.14 & TL\\
2 & 80 & 4, 7, 8 & 0.00 & 9.06 & 13.82 & TL & 5.46 & TL\\
2 & 70 & 5, 10, 2 & 0.00 & 2,006.56 & 14.38 & TL & 3.35 & TL\\
2 & 80 & 6, 9, 7 & 0.00 & 9.06 & 20.93 & TL & 5.66 & TL\\
2 & 70 & 7, 8, 1 & 0.00 & 17.11 & 15.64 & TL & 6.11 & TL\\
3 & 90 & 8, 9, 7 & 0.00 & 9.64 & 100.00 & TL & 100.00 & TL\\
 \midrule 
 \multicolumn{3}{r|}{\text{\diameter}} & 0.00 & 210.26 & 23.62 & 9,817.74 & 14.51 & TL \\ 
 \midrule 
3 & 100 & 1, 10, 6, 8 & 0.00 & 130.14 & 100.00 & TL & 100.00 & TL\\
2 & 80 & 1, 4, 6, 9 & 0.00 & 9.24 & 44.25 & TL & 100.00 & TL\\
3 & 80 & 2, 7, 3, 9 & 0.00 & 8.74 & 100.00 & TL & 100.00 & TL\\
3 & 90 & 4, 7, 8, 3 & 0.00 & 9.44 & 50.16 & TL & 100.00 & TL\\
 \midrule 
 \multicolumn{3}{r|}{\text{\diameter}} & 0.00 & 39.39 & 73.60 & TL & 100.00 & TL \\ 
\bottomrule
\end{tabular}
\label{tab:base:2dxintegrated_separate_gap_all}
\end{table}

In Table~\ref{tab:base:2dxintegrated_separate_gap_all}, we present the average percentage gap at the time limit of three hours and CPU times ($t$) in seconds across all technologies, when solving the lines in a combined fashion. We mark instances with ``TL'' when the time limit of three hours is reached.

For DBs, the results are exceptionally strong: all instances are solved to optimality and CPU times remain very low across all problem sizes, confirming the computational robustness of diesel fleets.
In contrast, BEBs show significant variability. While some small and medium-sized instances reach optimality or have low gaps, some medium and large cases exhibit very high gaps up to not even finding a feasible solution. The average optimality gap increases notably with instance size. CPU times for BEBs are consistently near the time limit, indicating that many instances are terminated prematurely, leaving room for improvement in both solution quality and computational performance.

FCEBs perform similar to BEBs, although for FCEBs, the gap size stays more consistent up until instances with 80 service trips, resulting in a smaller average optimality gap than for BEBs. 
For the largest instances, however, the FCEB gap also reaches 100\%, mirroring the BEB performance in the most complex cases.
Overall, this table reinforces that DBs remain computationally the most tractable, while FCEBs offer a promising balance in solvability at least for small and medium-sized instances, while BEBs pose significant computational challenges, particularly in larger integrated scheduling scenarios.

\subsection{Cold temperatures}
\label{rl:sec:comp:cold}

In this scenario, we assess the impact of cold temperatures on energy consumption across FCEBs and BEBs. The increase in energy consumption for BEBs is mainly based on empirical findings reported by \cite{ViriCity2020}, which analyzes the performance of 100 BEBs in varying temperature ranges. Specifically, during days with temperatures between –10°C and –15°C, the study observed an average increase in energy consumption of 14\% for 12-meter buses and 21\% for 18-meter buses. Given we assume that our fleet consists of 12-meter buses, we apply a 14\% increase in energy consumption to reflect the impact of cold temperatures. To do so, we set the energy consumption rate $\theta^{BEB}_{cold} = 1.14\theta$.

Although the study of \cite{ViriCity2020} only investigates BEBs, cold temperatures also negatively affect the energy efficiency of other propulsion systems. \cite{reithuber2025energy} investigate the impact of heating or cooling operations on energy and hydrogen consumption on FCEBs. They identify an increase of 8.4\% in total hydrogen consumption compared to the same driving cycle performed without heating. The hydrogen consumption is further increased by frequent door opening and results in an increase of 6.2\% compared to the driving cycle without door opening during cold temperatures. We therefore set the hydrogen consumption rate $\theta^{FCEB}_{cold} = 1.084 \cdot 1.062 \theta$.

\begin{table}[h!]
\centering
\scriptsize
\caption{Results for the original and the cold temperature scenario, solved separately and combined, grouped by technology.}
\begin{adjustbox}{max width=\textwidth}
\begin{tabular}{c r l| *{4}{c}| *{4}{c}| *{4}{c}| *{4}{c}}
\toprule
& & 
  & \multicolumn{8}{c|}{standard setting} 
  & \multicolumn{8}{c}{cold scenario}
\\ 

& & 
  & \multicolumn{4}{c|}{BEB} 
  & \multicolumn{4}{c|}{FCEB}
  & \multicolumn{4}{c|}{BEB}
 & \multicolumn{4}{c}{FCEB}\\
&  &  
  & \multicolumn{2}{c}{ separate} & \multicolumn{2}{c|}{ combined}
  & \multicolumn{2}{c}{ separate} & \multicolumn{2}{c|}{ combined}
  & \multicolumn{2}{c}{ separate} & \multicolumn{2}{c|}{ combined}
  & \multicolumn{2}{c}{ separate} & \multicolumn{2}{c}{ combined} \\
$|K|$ & $|V^I|$ & lines
   & \#b &  \#c &  \#b &  \#c 
      & \#b &  \#c &  \#b &  \#c 
   & \#b &  \#c &  \#b &  \#c 
  &  \#b &  \#c &  \#b &  \#c \\
\midrule
1 & 50 & 1, 10 & 13 & 18 & 12 & 18 & 8 & 17 & \textit{8} & \textit{17} & 17 & 24 & 14 & 25 & 10 & 20 & 8 & 20\\
1 & 30 & 1, 4 & 7 & 4 & 7 & 3 & 5 & 6 & 5 & 6 & 10 & 4 & 8 & 7 & 6 & 8 & 5 & 7\\
2 & 40 & 2, 7 & 7 & 6 & \textit{5} & \textit{8} & 6 & 8 & \textit{5} & \textit{7} & 7 & 9 & 5 & 13 & 6 & 9 & \textit{5} & \textit{8}\\
2 & 40 & 2, 9 & 8 & 12 & 7 & 10 & 7 & 9 & \textit{7} & \textit{8} & 9 & 12 & 9 & 10 & 7 & 10 & \textit{7} & \textit{10}\\
1 & 40 & 3, 9 & 8 & 13 & 8 & 10 & 7 & 9 & \textit{7} & \textit{9} & 9 & 13 & 9 & 12 & 7 & 11 & \textit{7} & \textit{10}\\
1 & 50 & 4, 7 & 9 & 7 & \textit{7} & \textit{10} & 7 & 9 & \textit{6} & \textit{8} & 10 & 9 & 8 & 15 & 7 & 11 & \textit{6} & \textit{10}\\
2 & 60 & 5, 10 & 14 & 20 & \textit{12} & \textit{20} & 9 & 19 & 8 & 18 & 16 & 28 & 13 & 31 & 10 & 21 & \textit{9} & \textit{20}\\
1 & 50 & 6, 9 & 10 & 12 & \textit{8} & \textit{10} & 9 & 11 & \textit{7} & \textit{10} & 11 & 13 & 9 & 13 & 9 & 12 & \textit{7} & \textit{11}\\
2 & 60 & 7, 8 & 9 & 7 & \textit{8} & \textit{8} & 8 & 10 & \textit{8} & \textit{9} & 9 & 10 & \textit{8} & \textit{9} & 8 & 11 & \textit{8} & \textit{11}\\
2 & 60 & 8, 9 & 10 & 13 & \textit{7} & \textit{11} & 9 & 11 & \textit{6} & \textit{10} & 11 & 13 & \textit{8} & \textit{14} & 9 & 12 & \textit{6} & \textit{13}\\
 \midrule 
 \multicolumn{3}{r|}{\text{\diameter}} & 9.5 & 11.2 & 8.1 & 10.8 & 7.5 & 10.9 & 6.7 & 10.2 & 10.9 & 13.5 & 9.1 & 14.9 & 7.9 & 12.5 & 6.8 & 12.0 \\ 
 \midrule 
2 & 70 & 1, 10, 6 & 17 & 19 & \textit{11} & \textit{21} & 12 & 21 & \textit{9} & \textit{20} & 21 & 27 & \textit{11} & \textit{30} & 14 & 24 & \textit{9} & \textit{22}\\
1 & 70 & 1, 4, 10 & 17 & 20 & 15 & 20 & 11 & 20 & \textit{11} & \textit{19} & 22 & 26 & 18 & 26 & 13 & 24 & \textit{11} & \textit{23}\\
3 & 70 & 2, 7, 9 & 13 & 17 & \textit{14} & \textit{15} & 11 & 15 & \textit{10} & \textit{15} & 14 & 19 & \textit{14} & \textit{19} & 11 & 17 & \textit{10} & \textit{16}\\
2 & 70 & 2, 9, 8 & 12 & 14 & \textit{10} & \textit{11} & 11 & 13 & \textit{8} & \textit{13} & 13 & 15 & \textit{14} & \textit{19} & 11 & 14 & \textit{8} & \textit{15}\\
2 & 60 & 3, 9, 5 & 12 & 17 & \textit{10} & \textit{14} & 10 & 14 & \textit{10} & \textit{13} & 13 & 19 & \textit{12} & \textit{18} & 10 & 16 & \textit{10} & \textit{15}\\
2 & 80 & 4, 7, 8 & 13 & 9 & \textit{10} & \textit{11} & 11 & 13 & \textit{9} & \textit{13} & 14 & 12 & \textit{10} & \textit{15} & 11 & 15 & \textit{9} & \textit{15}\\
2 & 70 & 5, 10, 2 & 16 & 21 & \textit{14} & \textit{25} & 11 & 21 & \textit{10} & \textit{19} & 18 & 30 & 15 & 27 & 12 & 23 & \textit{10} & \textit{22}\\
2 & 80 & 6, 9, 7 & 15 & 17 & \textit{12} & \textit{16} & 13 & 17 & \textit{10} & \textit{15} & 16 & 20 & \textit{13} & \textit{21} & 13 & 19 & \textit{10} & \textit{17}\\
2 & 70 & 7, 8, 1 & 12 & 9 & \textit{9} & \textit{12} & 10 & 13 & \textit{8} & \textit{13} & 14 & 12 & \textit{10} & \textit{15} & 11 & 15 & \textit{8} & \textit{14}\\
3 & 90 & 8, 9, 7 & 15 & 18 & - & - & 13 & 17 & - & - & 16 & 20 & \textit{20} & \textit{24} & 13 & 19 & - & -\\
 \midrule 
 \multicolumn{3}{r|}{\text{\diameter}} & 14.2 & 16.1 & 11.7 & 16.1 & 11.3 & 16.4 & 9.4 & 15.6 & 16.1 & 20.0 & 13.7 & 21.4 & 11.9 & 18.6 & 9.4 & 17.7 \\ 
 \midrule 
3 & 100 & 1, 10, 6, 8 & 21 & 21 & - & - & 16 & 25 & - & - & 25 & 30 & - & - & 18 & 28 & - & -\\
2 & 80 & 1, 4, 6, 9 & 17 & 16 & \textit{19} & \textit{18} & 14 & 17 & - & - & 21 & 17 & \textit{17} & \textit{22} & 15 & 20 & \textit{11} & \textit{19}\\
3 & 80 & 2, 7, 3, 9 & 15 & 19 & - & - & 13 & 17 & - & - & 16 & 22 & - & - & 13 & 20 & \textit{12} & \textit{20}\\
3 & 90 & 4, 7, 8, 3 & 15 & 11 & \textit{17} & \textit{26} & 13 & 15 & - & - & 16 & 15 & - & - & 13 & 18 & - & -\\
 \midrule 
 \multicolumn{3}{r|}{\text{\diameter}} & 17.0 & 16.8 & 18.0 & 22.0 & 14.0 & 18.5 & - & - & 19.5 & 21.0 & 17.0 & 22.0 & 14.8 & 21.5 & 11.5 & 19.5 \\ 

\bottomrule
\end{tabular}
\end{adjustbox}
\label{tab:cold:2dxintegrated_separateAll}
\end{table}

Table \ref{tab:cold:2dxintegrated_separateAll} presents the average total number of buses (\#b) and charging events (\#c) for BEBs, and FCEBs, under cold temperatures and compares them with the results from the standard setting. Again we look at two solution approaches: separately (each line group solved on its own) and combined (all lines in the group solved jointly). The results capture the impact of reduced energy efficiency in cold weather.
Each row corresponds to a specific instance defined by the number of depots ($|K|$), number of service trips ($|V^I|$), and the pair of lines considered. All instances not solved to optimality are displayed in italic. For instances, where we were not able to find a feasible solution, ``-'' is displayed.

In nearly all cases and over both technologies, solving lines combined reduces the number of buses compared to solving separately. This confirms the benefit of flexibility from joint planning, even under challenging conditions.

Identical to the standard setting, BEBs require the largest fleet size and the most recharging events on average across all combinations of lines. Cold temperatures reflect in increased battery demand and reduced battery performance, which results in a higher charging frequency. In some combined cases, where coordination becomes computationally heavy, no feasible solutions could be obtained. 
Compared to the standard setting, the required BEB fleet size remains unchanged in some instances but increases in the majority of cases under cold temperature conditions.

Even though FCEBs are much more resilient to cold than BEBs, FCEBs also experience some increase in the number of charging events compared to original conditions, but the increase is moderate and manageable: While FCEBs charge around 1.6 times on average in moderate temperatures, an average FCEB recharges at least 1.76 times when it is could outside. Nevertheless, the fleet size is almost identical across both scenarios. This demonstrates their suitability for environments with challenging weather.

\subsection{Battery preservation} 
\label{rl:sec:comp:battery}

To account for battery health and ensure realistic operational conditions, we constrain the minimum SOC of BEBs to remain between 30\% and 90\% of the maximum capacity throughout operation. Thus, we set the minimum SOC $s^{min}$ to $0.3s^{max}$, and the maximum SOC $s^{max}$ to 90\% of its original value. 

As shown by \cite{schmalstieg2014}, operating lithium-ion batteries at extreme SOC levels, particularly near full charge (above 90\%) or deep discharge (below 20–30\%), accelerates capacity fade and overall degradation. \cite{tomaszewska2019} further reinforce this point in their comprehensive review, highlighting how battery degradation accelerates significantly at both ends of the SOC spectrum, especially under high charging currents or frequent fast charging, which is common in electric bus systems.

Zheng et al. (2022) adopt the same 30–90\% SOC window in their en-route fast-charging optimization model, citing degradation-aware scheduling as a critical factor for extending battery lifetime and reducing long-term operational costs. They emphasize that maintaining SOC within this band helps simulate a quasi-linear charging region, which improves modeling accuracy and reflects real-world charging strategies where buses avoid prolonged charging near full capacity. 
\cite{He2020} similarly impose upper SOC limits in their fast-charging scheduling framework, aiming to balance charging efficiency and operational feasibility.

In addition to these technical motivations, our goal is to explore whether enforcing these SOC boundaries has a measurable effect on operational planning -- specifically, whether reduced battery availability leads to an increase in the required fleet size. By constraining the usable SOC range, we limit the energy buffer available for vehicle scheduling, which may in turn affect routing flexibility, charging frequency, and ultimately the number of buses needed to maintain service levels. Incorporating this constraint allows us to assess the trade-off between preserving battery health and maintaining operational efficiency in electric bus systems.

We now investigate whether these SOC constraints impact operational planning, specifically, whether reduced battery availability leads to a larger required fleet size due to limited scheduling flexibility and possible more frequent charging.

\begin{table}[ht]
\centering
\scriptsize
\caption{Results for the original and battery preservation scenario for BEBs, solved separately and combined.}
    \begin{tabular}{ccl|cccc|ccccrr}
\toprule
&&   & \multicolumn{4}{c|}{standard setting} 
  & \multicolumn{6}{c}{battery preservation scenario} \\
 &&& \multicolumn{2}{c}{separate} & \multicolumn{2}{c|}{combined} &
 \multicolumn{2}{c}{separate} & \multicolumn{2}{c}{combined}\\
$|K|$ & $|V^I|$ & lines
  &  \#b &  \#c &  \#b &  \#c & \#b &  \#c & \#b &  \#c & gap(\%) & $t$(s) \\
\midrule
1 & 50 & 1, 10 & 13 & 18 & 12 & 18 & 17 & 32 & 15 & 33 & 0.00 & 9.07\\
1 & 30 & 1, 4 & 7 & 4 & 7 & 3 & 12 & 8 & 10 & 10 & 0.00 & 10.25\\
2 & 40 & 2, 7 & 7 & 6 & \textit{5} & \textit{8} & 7 & 14 & 6 & 8 & 0.00 & 1,934.23\\
2 & 40 & 2, 9 & 8 & 12 & 7 & 10 & 10 & 16 & 9 & 11 & 0.00 & 77.64\\
1 & 40 & 3, 9 & 8 & 13 & 8 & 10 & 10 & 15 & 10 & 15 & 0.00 & 9.16\\
1 & 50 & 4, 7 & 9 & 7 & \textit{7} & \textit{10} & 12 & 14 & \textit{10} & \textit{17} & 0.17 & TL\\
2 & 60 & 5, 10 & 14 & 20 & \textit{12} & \textit{20} & 16 & 34 & 14 & 35 & 0.00 & 19.91\\
1 & 50 & 6, 9 & 10 & 12 & \textit{8} & \textit{10} & 12 & 16 & 10 & 16 & 0.00 & 170.74\\
2 & 60 & 7, 8 & 9 & 7 & \textit{8} & \textit{8} & 9 & 13 & \textit{8} & \textit{11} & 4.59 & TL\\
2 & 60 & 8, 9 & 10 & 13 & \textit{7} & \textit{11} & 12 & 15 & \textit{8} & \textit{14} & 15.50 & TL\\
 \midrule 
 \multicolumn{3}{r|}{\text{\diameter}} & 9.5 & 11.2 & 8.1 & 10.8 & 11.7 & 17.7 & 10.0 & 17.0 & 2.03 & 3,466.00 \\ 
 \midrule 
2 & 70 & 1, 10, 6 & 17 & 19 & \textit{11} & \textit{21} & 21 & 36 & \textit{12} & \textit{28} & 0.15 & TL\\
1 & 70 & 1, 4, 10 & 17 & 20 & 15 & 20 & 24 & 36 & 21 & 38 & 0.00 & 13.04\\
3 & 70 & 2, 7, 9 & 13 & 17 & \textit{14} & \textit{15} & 15 & 26 & \textit{14} & \textit{21} & 21.32 & TL\\
2 & 70 & 2, 9, 8 & 12 & 14 & \textit{10} & \textit{11} & 14 & 19 & \textit{9} & \textit{22} & 5.07 & TL\\
2 & 60 & 3, 9, 5 & 12 & 17 & \textit{10} & \textit{14} & 14 & 21 & 12 & 19 & 0.00 & 521.38\\
2 & 80 & 4, 7, 8 & 13 & 9 & \textit{10} & \textit{11} & 16 & 17 & \textit{10} & \textit{16} & 15.45 & TL\\
2 & 70 & 5, 10, 2 & 16 & 21 & \textit{14} & \textit{25} & 18 & 38 & 15 & 37 & 0.00 & 102.00\\
2 & 80 & 6, 9, 7 & 15 & 17 & \textit{12} & \textit{16} & 17 & 26 & \textit{13} & \textit{26} & 8.45 & TL\\
2 & 70 & 7, 8, 1 & 12 & 9 & \textit{9} & \textit{12} & 14 & 17 & \textit{9} & \textit{19} & 13.26 & TL\\
3 & 90 & 8, 9, 7 & 15 & 18 & - & - & 17 & 25 & - & - & 100.00 & TL\\
 \midrule 
 \multicolumn{3}{r|}{\text{\diameter}} & 14.2 & 16.1 & 11.7 & 16.1 & 17.0 & 26.1 & 12.8 & 25.1 & 16.37 & 7,633.96 \\ 
 \midrule 
3 & 100 & 1, 10, 6, 8 & 21 & 21 & - & - & 25 & 39 & \textit{18} & \textit{42} & 34.75 & TL\\
2 & 80 & 1, 4, 6, 9 & 17 & 16 & \textit{19} & \textit{18} & 24 & 24 & \textit{14} & \textit{22} & 13.57 & TL\\
3 & 80 & 2, 7, 3, 9 & 15 & 19 & - & - & 17 & 29 & \textit{18} & \textit{25} & 29.30 & TL\\
3 & 90 & 4, 7, 8, 3 & 15 & 11 & \textit{17} & \textit{26} & 18 & 20 & - & - & 100.00 & TL\\
 \midrule 
 \multicolumn{3}{r|}{\text{\diameter}} & 17.0 & 16.8 & 18.0 & 22.0 & 21.0 & 28.0 & 16.7 & 29.7 & 44.40 & TL \\ 

\bottomrule
\end{tabular}
\label{tab:battery:2dxintegrated_separate_all}
\end{table}

Table~\ref{tab:battery:2dxintegrated_separate_all} presents the results for BEBs under a battery preservation scenario, where charging operations are restricted to avoid deep charging, a typical strategy to extend battery life. It shows, for both the original and the battery preservation scenario, the average total number of BEBs required (\#b) and the corresponding number of charging events (\#c), comparing cases where lines are solved separately versus in combination. We mark instances with ``TL'' when the time limit of three hours is reached and with ``-'' when we could not find an integer feasible solution. Furthermore, all instances not solved to optimality are displayed in italic.

As expected, across all problem sizes, combining lines consistently reduces the number of BEBs needed, while the number of charging events tends to generally increase. This reflects the benefit of resource pooling: by jointly optimizing multiple lines, the scheduling algorithm can make more efficient use of the BEB fleet, reducing redundancy. 

Compared to the results for BEBs in the standard setting, the most notable difference is the increase in fleet size required in the battery preservation scenario. Across nearly all instance sizes, 
the number of BEBs needed is higher when battery preservation is prioritized. 
For combined instances with two and three lines, the average number of BEBs increases by 1.5.
This is a direct result of more conservative battery usage, which limits the allowable depth of discharge. Furthermore, a clear and important difference is also the increased number of charging events in the battery preservation setting. Due to the constraints imposed to protect battery health, buses need to recharge more frequently. This is reflected in the data: for nearly all instance sizes, the number of charging events in the battery preservation scenario exceeds that of the base scenario, often by a considerable margin.

In this scenario, we observe a significant improvement in computational performance compared to the base and cold scenarios. Feasible solutions are now found for instances with line 1, 10, 6, 8 and 2, 7, 3, 9, including the largest one. Moreover, optimality gaps decrease substantially across the board and the solver hits the time limit in a lot fewer cases for the battery preservation scenario compared to the original one. 
This demonstrates that the additional constraints, namely restricting charging to below 90\% and above 30\% of $s^{max}$, do not only promote battery health and longevity but also appear to simplify the underlying optimization problem. By narrowing the feasible region, these constraints reduce complexity and enable faster convergence, making large-scale electric fleet scheduling with battery-preserving policies much more tractable.
Overall, while the base scenario allows for tighter, more efficient scheduling and fleet reuse, the battery preservation scenario demands a larger fleet to maintain healthier charging behavior. This highlights a key operational trade-off between efficiency and long-term battery sustainability.

\section{Conclusion}
\label{conclusion}

In this work, we have modeled and solved the electric vehicle scheduling problem with multiple depots, multiple charging stations, and partial recharge. We have developed a graph representation that is an acyclic network and allows only time-feasible paths. Only two additional aspects need to be ensured for each vehicle schedule: the SOC of the vehicles along a path, ensuring that vehicles cannot run out of energy, and that each vehicle returns to the same depot as it started from.

While our 3-index MILP can be solved with any off-the-shelf solver directly, in order to accommodate multiple depots, constraints of exponential size are incorporated into our 2-index MILP formulation. We compare two types of these constraints and different tailored separation strategies. The best performing strategy for a lower number of depots relies on infeasible path constraints, separated only at new integer incumbent solutions during the execution of the branch-and-cut algorithm, while for a higher number of depots connectivity cuts, separated at new integer incumbent solutions with setting \texttt{All} works best. The 2-index-based branch-and-cut algorithm consistently solves more instances to optimality and in lower computation times than the 3-index model solved by CPLEX, for a low number of depots as well as for six and eight depots.

The comparison of realistic scenarios reveals important insights in the operational resilience and adaptability of each propulsion technology. 
Across all scenarios and technologies, we observe that combining multiple lines into a single scheduling problem reduces the total fleet size required. However, the number of recharging/refueling events decrease disproportionately on average, i.e., each individual vehicle charges more than in the separate case.

In general, serving the schedule requires at least as many, and often more, BEBs compared to DBs or FCEBs. BEBs are notably sensitive to environmental constraints, resulting in higher fleet requirements, whereas FCEBs can maintain stable fleet sizes even under cold temperatures. This indicates that FCEBs are more resilient to changes in operational conditions without significant performance degradation.

Interestingly, in the more restrictive cold temperature and battery preservation scenarios, despite requiring larger fleets and significantly more recharging, optimality gaps decrease for both BEBs and FCEBs. This suggests that tighter operational constraints can reduce the solution space and help guide the solver more efficiently toward high-quality solutions.

Taken together, these findings emphasize that integration between lines is a powerful lever to improve operational efficiency, but its implications, particularly regarding vehicle and infrastructure utilization, differ significantly between technologies.

Our approaches have been developed within the collaborative research project ZEMoS (Zero Emission Mobility Salzburg) and serve as decision support for fleet sizing decisions in two pilot regions in the country of Salzburg. 
Future work will involve the development of heuristic approaches for more complex problem versions, such as heterogeneous vehicle fleets.

\bibliographystyle{apalike}
\bibliography{Literature_Bus}
\end{document}